\documentclass[11pt]{article}
\usepackage{mathrsfs}
\usepackage{amsfonts}
\usepackage{bbm}
\usepackage[all]{xy}
\usepackage{eufrak}
\usepackage{amscd}
\usepackage{amsmath}
\def\O{{\cal O}}
\def\K{{\cal K}}

\input amssym.def
\input amssym.tex

\begin{document}

\begin{center}{\large\bf
Blocks with defect group $\mathbb Z_{2^n}\times \mathbb Z_{2^n}\times \mathbb Z_{2^m}$}

\bigskip\medskip{\large Chao Wu, Kun Zhang and Yuanyang Zhou}

\bigskip{\scriptsize School of Mathematics and
Statistics, Central China Normal University,
Wuhan, 430079, P.R. China
}
\end{center}

\bigskip\noindent{\bf Abstract} In this paper, we prove that a block algebra with defect group $\mathbb Z_{2^n}\times \mathbb Z_{2^n}\times \mathbb Z_{2^m}$, where $n\geq 2$ and $m$ is arbitrary, is Morita equivalent to its Brauer correspondent.

\bigskip\noindent{\bf Keywords}: Finite groups; blocks; Morita equivalences

\bigskip\bigskip\noindent{\bf 1. Introduction}

\bigskip Let $p$ be a prime number and $\O$ a complete discrete valuation ring with an algebraically closed field $k$ of characteristic $p$. Denote by $\K$ a fraction field of $\O$. Assume that $\K$ contains a $|G|$-th primitive root of unity for any finite group $G$ considered below when $\O$ has characteristic 0.

Let $G$ be a finite group, $b$ a block of $G$ over $\O$ with defect group $P$, and $c$ the Brauer correspondent of $b$ in the normalizer $N_G(P)$.
Assume that $P$ is abelian.
 {\it Brou\'e's abelian defect group conjecture} says that the block algebras $\O Gb$ and $\O N_G(P)c$ are derived equivalent (see \cite{B}). Denote by $\mathbb Z_n$ the residue class group modulo $n$. Recently, some authors proved the conjecture for some 2-blocks with abelian defect groups, such as $\mathbb Z_2\times \mathbb Z_2$ (see \cite{CEKL} and \cite{R}), $\mathbb Z_2\times \mathbb Z_2\times \mathbb Z_2$ (see \cite{E}), and $\mathbb Z_{2^n}\times \mathbb Z_{2^n}$ (see \cite{EKKS} for $n\geq 2$). The conjecture also has been verified in many other cases.

 We are especially interested in the case $\mathbb Z_{2^n}\times \mathbb Z_{2^n}$. In \cite{EKKS}, the authors classified 2-blocks of quasisimple groups with abelian defect groups, and then proved that the block algebras $\O Gb$ and $\O N_G(P)c$ are Morita equivalent when $P$ is $\mathbb Z_{2^n}\times \mathbb Z_{2^n}$ for some $n\geq 2$. By inspecting the classification, we observe that the block algebras $\O Gb$ and $\O N_G(P)c$ are Morita equivalent when $G$ is quasisimple and $P$ is $\mathbb Z_{2^n}\times\cdots \times \mathbb Z_{2^n}$ for some $n\geq 2$. So it may be interesting to ask whether such an observation holds for a general finite group $G$.
In the paper, we investigate blocks with defect group $\mathbb Z_{2^n}\times \mathbb Z_{2^n}\times \mathbb Z_{2^n}$, where $n\geq 2$. More generally, we prove Brou\'e's abelian defect group conjecture for blocks with defect groups $\mathbb Z_{2^n}\times \mathbb Z_{2^n}\times \mathbb Z_{2^m}$, where
$n\geq 2$ and $m$ is arbitrary.

\bigskip\noindent{\bf Theorem.}\quad {\it Assume that $P$ is $\mathbb Z_{2^n}\times \mathbb Z_{2^n}\times \mathbb Z_{2^m}$, where
$n\geq 2$ and $m$ is arbitrary. Then the block algebras $\O Gb$ and $\O N_G(P)c$ are Morita equivalent.}

\bigskip
Given a finite $p$-group $P$, {\it Donovan's conjecture} says that there are finitely many Morita equivalence classes of block algebras over $k$ with defect group $P$. It is expected that Donovan's conjecture should hold over $\O$ too. It is easy to see that the quotient group $N_G(P)/C_G(P)$ is trivial, or $\mathbb Z_3$, or $\mathbb Z_7$, or the Frobenius group $\mathbb Z_7\rtimes \mathbb Z_3$. By the structure theorem of block algebras with normal defect group, we conclude the following.

\bigskip\noindent{\bf Corollary.}\quad {\it Assume that $P$ is $\mathbb Z_{2^n}\times \mathbb Z_{2^n}\times \mathbb Z_{2^m}$, where
$n\geq 2$ and $m$ is arbitrary. Then the block algebra $\O Gb$ is Morita equivalent to $\O P$, or $\O (P\rtimes \mathbb Z_3)$, or $\O(P\rtimes \mathbb Z_7)$, or $\O(P\rtimes (\mathbb Z_7\rtimes \mathbb Z_3))$. In particular, Donovan's conjecture holds.  }

\bigskip
We remark that when $m\neq n$ the block algebra $\O Gb$ is Morita equivalent to $\O P$ or $\O (P\rtimes \mathbb Z_3)$, and that when $m=0$, by Proposition 6.7 below, the Morita equivalence between $\O Gb$ and $\O N_G(P)c$ in {\bf Theorem} is basic in the sense of \cite{P3} .

The proof of {\bf Theorem} heavily depends on $p$-extensions of certain perfect isometries (see Proposition 5.7 below). So the ring
$\O$ will be always assumed to have characteristic 0 in the rest of the paper except in the proof of {\bf Theorem}.

\bigskip\bigskip\noindent{\bf 2. Preliminaries}

\bigskip In this section, we collect some elementary lemmas.
Let $G$ and $G'$ be finite groups with common center $Z$. The central product $G\ast G'$ over $Z$ is isomorphic to the subgroup $G\otimes G'$ in the tensor product $\O G\otimes_{\O Z}\O G'$ and the group algebra $\O(G\ast G')$ is isomorphic to $\O G\otimes_{\O Z}\O G'$.

\bigskip\noindent{\bf Lemma 2.1.}\quad {\it Keep the notation as above. Then any block of $G\otimes G'$ is of the form $b\otimes b'$, where $b$ and $b'$ are blocks of $G$ and $G'$ respectively. If $P$ and $P'$ are defect groups of $b$ and $b'$ respectively, then $P\otimes P'$ is a defect group of $b\otimes b'$. Moreover, the order of the inertial quotient of the block $b\otimes b'$ is the product of the order of the inertial quotient of the block $b$ and that of $b'$.  }

\medskip\noindent{\it Proof.}\quad The map $G\times G'\rightarrow G\otimes G', (g, g')\mapsto g\otimes g'$ is a surjective group homomorphism with kernel $\{(z, z^{-1})|z\in Z\}$ isomorphic to $Z$. Denote by $Z'$ the maximal $p'$-subgroup of $Z$. Note that the quotient group $Z/Z'$ is a $p$-group. The homomorphism $G\times G'\rightarrow G\otimes G'$ factors through the canonical homomorphisms \begin{center}$G\times G'\rightarrow (G\times G')/Z'$ and $(G\times G')/Z'\rightarrow G\otimes G'$. \end{center}

Let $a\in \O G$ and $a'\in \O G'$. We denote by $a\tilde\otimes a'$ the element of the tensor product $\O G\otimes_\O \O G'$ determined by $a$ and $a'$, in order to differentiate $a\tilde\otimes a'$ and $a\otimes a'$. There is an $\O$-algebra isomorphism $\O(G\times G')\cong \O G\otimes_\O \O G'$ sending $(x_1, x_2)$ onto $x_1\tilde\otimes x_2$ for any $x\in G$ and any $x'\in G'$. We identify the algebras $\O(G\times G')$ and $\O G\otimes_\O \O G'$. Then any block of $G\times G'$ is of form $b\tilde\otimes b'$, where $b$ and $b'$ are blocks of $G$ and $G'$. If $P$ and $P'$ are defect groups of $b$ and $b'$, then $P\times P'$ is a defect group of $b\tilde\otimes b'$.
The inertial quotient of $b\tilde\otimes b'$ is the direct product of the inertial quotients of $b$ and $b'$.

Set $e=\frac{1}{|Z'|}\sum_{z\in Z'} z$. The homomorphism $G\times G'\rightarrow (G\times G')/Z'$ induces an algebra isomorphism $\O(G\times G')e\cong \O ((G\times G')/Z'). $ We identify the two algebras $\O(G\times G')e$ and $\O ((G\times G')/Z')$.
By the last paragraph, any block of $(G\times G')/Z'$ is of form $b\tilde\otimes b'$, where $b$ and $b'$ are blocks of $G$ and $G'$; if $P$ and $P'$ are defect groups of $b$ and $b'$, then $P\times P'$ is isomorphic to a defect group of $b\tilde\otimes b'$ (see \cite[Theorem 5.7.4]{NT}); the inertial quotient of $b\tilde\otimes b'$ is isomorphic to the direct product of the inertial quotients of $b$ and $b'$ (see \cite[Theorem 3.6]{P8}).

The homomorphism $(G\times G')/Z'\rightarrow G\otimes G'$ induces a surjective homomorphism $\O ((G\times G')/Z')\rightarrow \O (G\otimes G')$ with kernel the central $p$-subgroup $Z/Z'$, which induces a one to one correspondence between blocks of $(G\times G')/Z'$ and $G\otimes G'$ (see \cite[5.8.11]{NT}). If $D/Z'$ is a defect group of a block of $(G\times G')/Z'$, then $D/Z$ is a defect group of the corresponding block of $G\otimes G'$ (see \cite[Theorem 5.8.10]{NT}). By \cite[Theorem 2.9]{P7}, a block of $(G\times G')/Z'$ and the corresponding block of $G\otimes G'$ have isomorphic inertial quotients. The proof is done.

\medskip\noindent{\bf 2.2.}\quad Let ${\rm CF}(G, \K)$ be the $\K$-vector space of all $\K$-valued class functions on $G$, and ${\rm CF}_{p'}(G, \K)$ the subspace of all functions in ${\rm CF}(G, \K)$ vanishing on $p$-singular elements of $G$.
For any $\chi\in {\rm CF}(G, \K)$ and any central idempotent $e$ in $\O G$, extending $\chi$ to a function over $\K G$ by the $\K$-linearity, we define a new class function $e\cdot \chi$ on $G$ by the equality $(e\cdot \chi)(g)=\chi(g e)$ for any $g\in G$. We set ${\rm CF}(G, e, \K)=\{e\cdot \chi|\chi\in {\rm CF}(G, \K)\}$ and ${\rm CF}_{p'}(G, e, \K)={\rm CF}(G, e, \K)\cap {\rm CF}_{p'}(G, \K)$. It is easy to see that ${\rm CF}(G, e, \K)\,\mbox{\rm and}\, {\rm CF}_{p'}(G, e, \K)$ are subspaces of ${\rm CF}(G, \K)$ and that we have direct sum decompositions
$${\rm CF}(G, \K)=\bigoplus_b {\rm CF}(G, b, \K)\,\mbox{\rm and}\, {\rm CF}_{p'}(G, \K)=\bigoplus_b {\rm CF}_{p'}(G, b, \K),$$ where $b$ runs over the set of blocks of $G$. Denote by ${\rm Irr}(G)$ and ${\rm Irr}(G, b)$ the set of irreducible ordinary characters of the group $G$ and the set of irreducible ordinary characters in a block $b$ of $G$. It is known that ${\rm Irr}(G)$ and ${\rm Irr}(G, b)$ are bases of ${\rm CF}(G, \K)$ and ${\rm CF}(G, b, \K)$ and that the orthogonality relationship of characters determines the inner products on ${\rm CF}(G, \K)$ and ${\rm CF}(G, b, \K)$. We denote by ${\rm CF}_{p'}(G, \O)$ the $\O$-submodule of all $\O$-valued class functions on $G$ vanishing on all $p$-singular elements of $G$, and by ${\rm CF}_{p'}(G,\,e,\, \O)$ the intersection of ${\rm CF}_{p'}(G, \O)$ and ${\rm CF}(G, e, \K)$.

\medskip\noindent{\bf 2.3.}\quad Let $b$ and $b'$ be blocks of $G$ and $G'$. Denote by $b'^\circ$ the inverse image of $b'$ through the opposite isomorphism $\O G'\cong \O G'$. There is an $\O$-algebra isomorphism $\O(G\times G')\cong \O G\otimes_\O \O G'$. Identifying the two algebras $\O(G\times G')$ and $\O G\otimes_\O \O G'$, $b\otimes b'^\circ$ is a block of $G\times G'$. Given a generalized character $\nu$ of $G\times G'$ in the block $b\otimes b'^\circ$, we define a map $I_\nu: \mathbb Z {\rm Irr}(G', b')\rightarrow \mathbb Z{\rm Irr}(G, b)$ by the equality $$I_\nu(\alpha)(g)=\frac{1}{|G'|}\sum_{h\in G'}\nu(g, h)\alpha(h)$$ for any $\alpha\in \mathbb Z {\rm Irr}(G', b')$ and any $g\in G$, where $\mathbb Z{\rm Irr}(G, b)$ and $\mathbb Z{\rm Irr}(G', b')$ are the free abelian groups generated by ${\rm Irr}(G, b)$ and by ${\rm Irr}(G', b')$. The map $I_\nu$ can be extended to a $\K$-linear map ${\rm CF}(G', b', \K)\rightarrow {\rm CF}(G, b, \K)$, denoted by $I_\nu^\K$. We remind that the definition of $I_\nu$ is slightly different from the definition of $I_\mu$ in \cite[($\mathcal F$)]{B1}. We modify the definition of $I_\mu$ in \cite[($\mathcal F$)]{B1}, mainly to fill gaps, which are produced when we use to $\mu$ the induction and extension of generalized characters. Such a modification does not influence the use of all known theorems on perfect isometries (see \cite{B1}). If $I_\nu$ is a perfect isometry and we denote by $\chi^*$ the dual of an ordinary character $\chi$, then $$\nu=\sum_{\chi\in {\rm Irr}(G',\, b')} I_\nu(\chi)\times \chi^*,$$ where $I_\nu(\chi)\times \chi^*$ denotes the generalized ordinary character of $G\times G'$ determined by $\chi^*$ and $I_\nu(\chi)$.

\medskip\noindent{\bf 2.4.}\quad Let $u$ be a $p$-element of $G$. We denote by $d_G^u$ the surjective map ${\rm CF}(G, \K)\rightarrow {\rm CF}_{p'}(C_G(u), \K)$ defined by $d_G^u(\chi)(s)=\chi(us)$ for any $\chi\in {\rm CF}(G, \K)$ and any $p$-regular element $s$ of $C_G(u)$, and by $e_G^u$ the map ${\rm CF}_{p'}(C_G(u), \K)\rightarrow {\rm CF}(G, \K)$ sending any $\varphi$ onto the function $e_G^u(\varphi)$ which takes $\varphi(s)$ at $us$ for any $p$-regular element $s$ of $C_G(u)$ and $0$ at $g\in G$ if the $p$-part of $g$ is not conjugate to $u$. The composition $d_G^u\circ e_G^u$ is the identity map on ${\rm CF}_{p'}(C_G(u), \K)$. Let $(u,\,e)$ be a $b$-Brauer element. We denote by $d_G^{(u,\, e)}$ the composition of $d_G^u$ and the projection ${\rm CF}_{p'}(C_G(u), \K)\rightarrow {\rm CF}_{p'}(C_G(u), e, \K)$. By \cite[Theorem A2.1]{B2} $e_G^u$ maps ${\rm CF}_{p'}(C_G(u), e, \K)$ into ${\rm CF}(G, b, \K)$ and we denote by $e_G^{(u,\,e)}$ the restriction of $e_G^u$ on ${\rm CF}_{p'}(C_G(u), e, \K)$. The composition $d_G^{(u,\,e)}\circ e_G^{(u,\,e)}$ is the identity map on ${\rm CF}_{p'}(C_G(u), e, \K)$.

\medskip\noindent{\bf 2.5.}\quad Let $(P, f)$ be a maximal $b$-Brauer pair. For any subgroup $Q$ of $P$, there is a unique $b$-Brauer pair $(Q, f_Q)$ contained in $(P, f)$. We denote by ${\rm Br}(P, f)$ the Brauer category (see \cite{T}), whose objects are all $b$-Brauer pairs $(Q, f_Q)$ contained in $(P, f)$ and whose morphisms from $(Q, f_Q)$ to $(R, f_R)$ are group homomorphisms $Q\rightarrow R$ induced by elements $g\in G$ such that $(Q, f_Q)^g\leq (R, f_R)$. We assume that the blocks $b$ and $b'$ have a common defect group $P$ and that the inclusions $P\subset G$ and $P\subset G'$ induce an isomorphism between the categories ${\rm Br}(P, f)$ and ${\rm Br}(P, f')$, where $(P, f')$ is a suitable maximal $b'$-Brauer pair. The map $I_\nu^\K:{\rm CF}(G',\, b',\, \K)\rightarrow {\rm CF}(G,\, b,\, \K)$ is said to be {\it compatible with a local system} $\{I_T: {\rm CF}_{p'}(C_{G'}(T),\, f'_T,\, \K)\rightarrow {\rm CF}_{p'}(C_G(T),\, f_T, \,\K)\}_{\{T ({\rm cyclic})\subset P\}}$ (see \cite[Definition 4.3]{B1}), if for any cyclic subgroup $T$ of $P$ and any generator $x$ of $T$, we have $d_G^{(x,\, f_T)}\circ I_\nu^\K=I_T\circ d_{G'}^{(x,\, f'_T)}$.
By the isomorphism theorem of groups, the local system is uniquely determined by $I_\nu^\K$ if it exists.

\bigskip\noindent{\bf Lemma 2.6.}\quad {\it Keep the notation and the assumption of 2.5 and assume that $P$ is abelian and that the map $I_{\nu}^\K:{\rm CF}(G',\, b',\, \K)\rightarrow {\rm CF}(G,\, b,\, \K)$ is compatible with a local system $\{I_T\}_{\{T ({\rm cyclic})\subset P\}}$. Then we have $I_\nu^\K=\sum_{x\in J} e_G^{(x,\, f_{\langle x \rangle})}\circ I_{\langle x \rangle}\circ d_{G'}^{(x,\, f'_{\langle x \rangle})}$, where $\langle x \rangle$ is the cyclic subgroup generated by $x$ and $J$ is a set of representatives of orbits of the conjugation action of $N_G(P, f)$ on $P$. }

\medskip\noindent{\it Proof.}\quad Obviously $\{(x, f_{\langle x\rangle})|x\in J\}$ is a set of representatives of orbits of the action of $N_G(P, f)$ on the set of $b$-Brauer elements contained in $(P, f)$, and since $P$ is abelian, $\{(x, f_{\langle x \rangle})|x\in J\}$ is a set of representatives of $G$-conjugacy classes of $b$-Brauer elements. The family $\{d_{G}^{(x,\, f_{\langle x \rangle})}|x\in J\}$ induces a $\K$-linear isomorphism ${\rm CF}(G, b, \K)\cong \bigoplus_{x\in J}{\rm CF}_{p'}(C_G(x), f_{\langle x \rangle}, \K)$ (see the last paragraph of \cite[4A]{B1}).
Since $d_{G}^{(x',\, f_{\langle x' \rangle})}\circ e_{G}^{(x,\, f_{\langle x \rangle})}=0$ for different $x, x'\in J$, we have $$d_{G}^{(x,\, f_{\langle x \rangle})}\circ I_\nu^\K=d_{G}^{(x,\, f_{\langle x \rangle})}\circ \Big(\sum_{y\in J} e_G^{(y,\, f_{\langle y\rangle})}\circ I_{\langle y \rangle}\circ d_{G'}^{(y,\, f'_{\langle y \rangle})}\Big).$$ So $I_\nu^\K=\sum_{y\in J} e_G^{(y,\, f_{\langle y\rangle})}\circ I_{\langle y \rangle}\circ d_{G'}^{(y,\, f'_{\langle y \rangle})}$.

\bigskip The main idea of the following proposition is already included in the first sentence of \cite[1.6]{PU}.

\bigskip\noindent{\bf Proposition 2.7.}\quad {\it Keep the notation and the assumption of 2.5 and assume that $P$ is abelian. Then the map $I_{\nu}^\K: {\rm CF}(G',\, b',\, \K)\rightarrow {\rm CF}(G,\, b,\, \K)$ is compatible with a local system $\{I_T\}_{\{T ({\rm cyclic})\subset P\}}$ if and only if $I_\nu^\K(\lambda\ast\chi)=\lambda\ast I_\nu^\K(\chi)$ for any $\chi\in {\rm CF}(G', b', \K)$ and any $N_G(P,\,f)$-stable character $\lambda$ of $P$, where $\ast$ denotes the $\ast$-construction of characters due to Brou\'e and Puig. }

\bigskip Before the proof, we remark that any ordinary character $\lambda$ of $P$ is $N_G(P,\,f)$-stable if and only if it is $N_{G'}(P,\, f')$-stable, since the Brauer categories ${\rm Br}(P,\, f)$ and ${\rm Br}(P,\, f')$ are isomorphic.

\bigskip\noindent{\it Proof.}\quad
We assume that the map $I_{\nu}$ is compatible with a local system $\{I_T\}_{\{T ({\rm cyclic})\subset P\}}$. For any $b$-Brauer element $(u, f_{\langle u \rangle})$ and any $\chi\in {\rm CF}(G, b, \K)$, set $\chi^{(u,\, f_{\langle u \rangle})}=(e^{(u,\, f_{\langle u \rangle})}\circ d^{(u,\, f_{\langle u \rangle})})(\chi)$.
Let $J$ be a set of representatives of orbits of the conjugation action of $N_G(P, f)$ on $P$. Since $P$ is abelian, $\{(x, f_{\langle x \rangle})|x\in J\}$ is a set of representatives of $G$-conjugacy classes of $b$-Brauer elements. Since the categories ${\rm Br}(P, f)$ and ${\rm Br}(P, f')$ are isomorphic, $\{(x, f'_{\langle x \rangle})|x\in J\}$ is a set of representatives of orbits of the action of $N_{G'}(P, f')$ on the set of $b'$-Brauer elements contained in $(P, f')$ and $\{(x, f'_{\langle x \rangle})|x\in J\}$ is also a set of representatives of the $G'$-conjugacy classes of $b'$-Brauer elements.
For any $\chi\in {\rm CF}(G', b', \K)$, we have $\chi=\sum_{x\in J}\chi^{(x,\, f'_{\langle x \rangle})}$. Since
$I_\nu^\K(\lambda\ast\chi)^{(x,\, f_{\langle x \rangle})}  =   \lambda(x)e_{G}^{(x,\, f_{\langle x \rangle})}\circ I_{\langle x \rangle}\circ d_{G'}^{(x,\, f'_{\langle x \rangle})}(\chi^{(x,\, f'_{\langle x \rangle})})  =(\lambda\ast I_\nu^\K(\chi))^{(x,\, f_{\langle x \rangle})}$ (see Lemma 2.6), we have $$I_\nu^\K(\lambda\ast\chi)=\sum_{x\in J} I_\nu^\K(\lambda\ast\chi)^{(x,\, f_{\langle x \rangle})}=\sum_{x\in J}(\lambda\ast I_\nu^\K(\chi))^{(x,\, f_{\langle x \rangle})}=\lambda\ast I_\nu^\K(\chi).$$

Conversely, we assume that $I_\nu^\K(\lambda\ast\chi)=\lambda\ast I_\nu^\K(\chi)$ for any $\chi\in {\rm CF}(G', b', \K)$ and any $N_G(P,\, f)$-stable character $\lambda$ of $P$.
Obviously $N_G(P, f)$ acts on ${\rm Irr}(P)$ and by Brauer's permutation Lemma, the number of orbits of the action of $N_G(P, f)$ on ${\rm Irr}(P)$ is equal to the cardinality of $J$. We list the orbit sums of the action of $N_G(P, f)$ on ${\rm Irr}(P)$ as $\lambda_j$, where $1\leq j\leq n$ and $ n=|J|$. For any $\lambda_i$ and any $\chi\in {\rm Irr}(G', b')$, we have $\lambda_i\ast\chi=\sum_{u\in J} \lambda_i(u)\chi^{(u,\, f'_{\langle u \rangle})}.$
The matrix $(\lambda_i(u))$ is an invertible matrix, since ${\rm Irr}(P)$ is a basis of ${\rm CF}(P, \K)$. Denoting by $(a_{iu})$ the inverse of $(\lambda_i(u))$, we have $$\chi^{(u,\, f'_{\langle u \rangle})}=\sum_i a_{iu}\lambda_i\ast\chi. \leqno 2.7.1$$
Since $I_\nu^\K(\lambda_i\ast\chi)=\lambda_i\ast I_\nu^\K(\chi)=\sum_{u\in J} \lambda_i(u)(I_\nu^\K(\chi))^{(u,\, f_{\langle u \rangle})}$, similarly we have $$(I_\nu^\K(\chi))^{(u,\, f_{\langle u \rangle})}=\sum_i a_{iu}I_\nu^\K(\lambda_i\ast\chi).\leqno 2.7.2$$ Combining the equalities 2.7.1 and 2.7.2, we have $$(I_\nu^\K(\chi))^{(u,\, f_{\langle u \rangle})}=I_\nu^\K(\chi^{(u,\, f'_{\langle u \rangle})}).\leqno 2.7.3$$

For any $u\in J$, we denote by $V^{(u,\, f_{\langle u \rangle})}$ the image of the homomorphism $e^{(u,\, f_{\langle u \rangle})}$. Obviously the set $\{\chi^{(u,\, f_{\langle u \rangle})}| \chi\in {\rm Irr}(G, b)\}$ is a generator set of the space $V^{(u,\, f_{\langle u \rangle})}$, the space ${\rm CF}(G, b, \K)$ is equal to the sum $\bigoplus_{u\in J} V^{(u,\, f_{\langle u \rangle})}$
since the family $\{d_{G}^{(u,\, f_{\langle u \rangle})}|u\in J\}$ induces a $\K$-linear isomorphism ${\rm CF}(G, b, \K)\cong \bigoplus_{u\in J}{\rm CF}_{p'}(C_G(u), f_{\langle u \rangle}, \K)$, and the kernel of $d_{G}^{(u,\, f_{\langle u \rangle})}$ is equal to the sum $\bigoplus_{v\in J-\{u\}} V^{(v,\, f_{\langle v\rangle})}$.
By the equality 2.7.3, $I_\nu^\K$ maps $V^{(u,\, f'_{\langle u \rangle})}$ onto $V^{(u,\, f_{\langle u \rangle})}$ and thus it maps the kernel of $d_{G'}^{(u,\, f'_{\langle u \rangle})}$ onto the kernel of $d_{G}^{(u,\, f_{\langle u \rangle})}$. Therefore $I_\nu^\K$ induces a $\K$-linear homomorphism $$I_{\langle u \rangle}: {\rm CF}_{p'}(C_{G'}(u), f'_{\langle u \rangle}, \K)\rightarrow {\rm CF}_{p'}(C_G(u), f_{\langle u \rangle}, \K)$$ such that $d_{G}^{(u,\, f_{\langle u \rangle})}\circ I_\nu^\K=I_{\langle u \rangle}\circ d_{G'}^{(u,\, f'_{\langle u \rangle})}$.

\bigskip\noindent{\bf Lemma 2.8.}\quad {\it Keep the notation and the assumption of 2.5 and assume that $P$ is abelian. Assume that the map $I_\nu^\K:{\rm CF}(G',\, b',\, \K)\rightarrow {\rm CF}(G,\, b,\, \K)$ is compatible with a local system $\{I_T\}_{\{T ({\rm cyclic})\subset P\}}$ and that $I_\nu^\K$ maps ${\rm CF}(G',\, b',\, \O)$ into ${\rm CF}(G,\, b,\, \O)$. Then for any cyclic subgroup $T$ of $P$, $I_T$ maps ${\rm CF}_{p'}(C_{G'}(T),\, f'_T,\, \O)$ into ${\rm CF}_{p'}(C_G(T),\, f_T, \,\O)$. }

\medskip\noindent{\it Proof.}\quad Take $x\in P$. Since $e_G^{(x,\, f_{\langle x \rangle})}\circ I_{\langle x \rangle}\circ d_{G'}^{(x,\, f_{\langle x \rangle}')}=e_G^{(x,\, f_{\langle x \rangle})}\circ d_G^{(x,\, f_{\langle x \rangle})}\circ I_\nu^\K=I_\nu^\K\circ e_{G'}^{(x,\, f_{\langle x \rangle}')}\circ d_{G'}^{(x,\, f_{\langle x \rangle}')}$ (the second equality is obtained by the equality 2.7.3), we have $e_G^{(x,\, f_{\langle x \rangle})}\circ I_{\langle x \rangle}=I_\nu^\K\circ e_{G'}^{(x,\, f_{\langle x \rangle}')}$. Since the homomorphisms $e^{(x, \,f_{\langle x \rangle})}_G$ maps ${\rm CF}_{p'}(C_G({\langle x \rangle}),\, f_{\langle x \rangle},\, \O)$ into ${\rm CF}(G,\, b,\, \O)$ and $I_\nu^\K$ maps ${\rm CF}(G',\, b',\, \O)$ onto ${\rm CF}(G,\, b,\, \O)$, the image of ${\rm CF}_{p'}(C_{G'}({\langle x \rangle}),\, f'_{\langle x \rangle},\, \O)$ through $e_{G}^{(x,\, f_{\langle x \rangle})}\circ I_{\langle x \rangle}$ is contained in ${\rm CF}(G,\, b,\, \O)$. Since $$d_{G}^{(x,\, f_{\langle x \rangle})}({\rm CF}(G,\, b,\, \O))\subset {\rm CF}_{p'}(C_{G}({\langle x \rangle}),\, f_{\langle x \rangle},\, \O),$$ $I_{\langle x \rangle}$ maps ${\rm CF}_{p'}(C_{G'}({\langle x \rangle}),\, f'_{\langle x \rangle},\, \O)$ into ${\rm CF}_{p'}(C_G({\langle x \rangle}),\, b_{\langle x \rangle},\, \O)$.

\bigskip\noindent{\bf 2.9.}\quad  Now we take an $(\O G,\, \O G')$-bimodule $M$ inducing a Morita equivalence between $\O Gb$ and $\O G'b'$. By \cite[Theorem 2.1 and Exercise 4.5]{AF}, the Morita equivalence induced by $M$ induces an algebra isomorphism $$\rho:Z(\O G'b')\cong Z(\O Gb)$$ such that for any $a\in Z(\O Gb)$ and any $a'\in Z(\O G'b')$, $a$ and $a'$ correspond to each other if and only if $am=ma'$ for any $m\in M$, where $Z(\O Gb)$ and $Z(\O G'b')$ are the centers of the algebras $\O Gb$ and $\O'G'b'$. In order to consider the character of the bimodule $M$ and basic Morita equivalences in Section 3, we also regard the module $M$ as an $\O( G\times G')$-bimodule by the equality $(g, g'^{-1})m=gmg'$ for any $g\in G$, any $g'\in G'$ and any $m\in M$.
Let $\mu$ be the character of $M$. The map $I_{\mu}: \mathbb Z{\rm Irr}(G', b')\rightarrow \mathbb Z{\rm Irr}(G, b)$ is a perfect isometry. By \cite[Theorem 1.5]{B1}, the perfect isometry $I_\mu$ induces an algebra isomorphism $$\rho':Z(\O G'b')\cong Z(\O Gb).$$

\bigskip\noindent{\bf Lemma 2.10.}\quad {\it Keep the notation above. The isomorphisms $\rho$ and $\rho'$ coincide. }

\medskip\noindent{\it Proof.}\quad Consider the extensions of the isomorphisms $\rho$ and $\rho'$ by the $\K$-linearity \begin{center}$\rho^\K: Z(\K G' b')\cong Z(\K Gb)$ and $\rho'^\K: Z(\K G' b')\cong Z(\K Gb)$. \end{center} For any $\chi\in {\rm Irr}(G)$, we denote by $e_\chi$ the central primitive idempotent $\frac{\chi(1)}{|G|}\sum_{g\in G}\chi(g)g^{-1}$ in $\K G$.  Set $e_\chi=\rho^\K(e_{\chi'})$ for any $\chi'\in {\rm Irr}(G',\,b')$. We have $e_\chi m=m e_{\chi'}$ for any $m\in M$, and thus $\mu(e_\chi,\, 1)=\mu(1,\, e_{\chi'^*})\neq 0$. Since $\mu=\sum_{\chi'\in {\rm Irr}(G',\, b')}(I_\mu(\chi')\times \chi'^*)$, $I_\mu(\chi')(e_{\chi})$ is nonzero and $I_\mu(\chi')$ is equal to $\chi$.
On the other hand, by the proof of \cite[Theorem 1.5]{B1} the isomorphism $\rho'^\K$ also sends $e_{\chi'}$ onto $e_{\chi}$. So $\rho'^\K=\rho^\K$.

\bigskip\bigskip\noindent{\bf 3. Morita equivalences and the $\ast$-structure}

\bigskip In this section, we will prove that if two block algebras with abelian defect groups are basically Morita equivalent, then there is a Morita equivalence between the two corresponding block algebras inducing a perfect isometry compatible with a local system (see Proposition 3.6).
In order to do that, we firstly recall the computation of generalized decomposition numbers in \cite{P0}.

\medskip\noindent{\bf 3.1.}\quad Throughout the section, all $\O$-algebras and $\O$-modules are $\O$-free of finite $\O$-rank. For any ring $R$, denote by $R^*$ the multiplicative group of $R$. Let $G$ be a finite group. An $\O$-algebra $A$ is an {\it interior $G$-algebra} if there is a group homomorphism $G\rightarrow A^*$. The conjugation of $G$ in $A$ induces a group homomorphism $G\rightarrow {\rm Aut}(A)$ so that $A$ becomes a $G$-algebra, where ${\rm Aut}(A)$ denotes the automorphism group of the $\O$-algebra $A$.
For any subgroup $H$ of $G$, we denote by $A^H$ the subalgebra of all $H$-fixed elements in the $G$-algebra $A$. A {\it pointed group} on $A$ is a pair $(H, \alpha)$, where $H$ is a subgroup of $G$ and $\alpha$ is an $(A^H)^*$-conjugacy class of primitive idempotents in $A^H$. We often write $(H, \alpha)$ as $H_\alpha$, and say that $\alpha$ is a {\it point} of $H$ on $A$. The pointed group $H_\alpha$ is {\it contained} in another pointed group $K_\beta$, denoted by $H_\alpha\leq K_\beta$, if $H$ is a subgroup of $K$ and there are $i\in \alpha$ and $j\in \beta$ such that $ji=ij=i$.

\medskip\noindent{\bf 3.2.}\quad Let $H$ and $K$ be two subgroups of $G$ such that $K\leq H$. We denote by ${\rm Tr}^H_K: A^K\rightarrow A^H$ the relative trace map and by $A^H_K$ the image of ${\rm Tr}^H_K$. We consider the quotient $$A(H)=A^H\Big/\Big(J(\O)A^H+\sum_K A_K^H  \Big)$$ where $K$ runs over all proper subgroups of $H$, and denote by ${\rm Br}_H^A$ the canonical surjective algebra homomorphism from $A^H$ to $A(H)$. If $A(H)$ is not zero, then $H$ has to be a $p$-subgroup. Clearly the inclusion $G\subset \O G$ endows $\O G$ an interior $G$-algebra structure. Let $Q$ be a $p$-subgroup of $G$. Then the inclusion $\O C_G(Q)\subset (\O G)^Q$ induces an algebra isomorphism $(\O G)(Q)\cong kC_G(Q)$. We identify $(\O G)(Q)$ and $k C_G(Q)$ in the sequel.

\medskip\noindent{\bf 3.3.}\quad
A pointed group $H_\alpha$ on $A$ is {\it local} if $H$ is a $p$-subgroup of $G$ and the image ${\rm Br}_H^A(\alpha)\neq \{0\}$. In this case, $\alpha$ is called a local point of $H$ on $A$, and since ${\rm Br}_H^A$ is a surjective algebra homomorphism, ${\rm Br}_H^A(\alpha)$ is a set of primitive idempotents in $A(H)$ and by the lifting theorem of idempotents, the correspondence $\alpha\mapsto {\rm Br}_H^A(\alpha)$ gives a bijection between the local points of $H$ on $A$ and the points of $H$ on $A(H)$. Let $x$ be an element of $G$. We call $x_\alpha$ a local pointed element if $\langle x\rangle_\alpha$ is a local pointed group on $A$. The {\it local pointed element} $x_\alpha$ is contained in a pointed group $K_\beta$ if $\langle x\rangle_\alpha$ is contained in $K_\beta$.

\medskip\noindent{\bf 3.4.}\quad Now we consider the group algebra $\O G$. Let $b$ a block of $G$. Then $\{b\}$ is a point of $G$ on $\O G$.
Let $x_\alpha$ be a local pointed element contained in the pointed group $G_{\{b\}}$ on $\O G$, and take $j\in \alpha$. Let $\chi\in {\rm Irr}(G, b)$, afforded by a $\K G$-module $M$. The product $jM$ of $j$ and $M$ is an $\O \langle u\rangle$-module, whose character is denoted by $\chi^\alpha$.
Denote by ${\cal P}_{\O G}(x)$ the set of all $\alpha$ such that $x_\alpha$ is a local pointed element on $\O G$. By \cite[Corollary 4.4]{P0}, for any $p'$-element $s$ of $C_G(x)$, we have $$\chi(xs)=\sum_{\alpha\in {\cal P}_{\O G}(x)}\chi^\alpha(x)\varphi_\alpha(s). \leqno 3.4.1$$ If $x_\alpha$ is not contained in $G_{\{b\}}$, then $\chi^\alpha(x)$ is equal to 0.

\medskip\noindent{\bf 3.5.}\quad
By \cite[Theorem 1.2]{P0}, all maximal local pointed groups contained in $G_{\{b\}}$ are $G$-conjugate to each other. Let $P_\gamma$ be a maximal local pointed group contained in $G_{\{b\}}$, choose $i\in \gamma$ and set $A_\gamma=i Ai$. Obviously $A_\gamma$ is an $\O$-algebra and admits a group homomorphism $\psi:P\rightarrow A_\gamma^*, u\mapsto ui$. So $A_\gamma$ is an interior $P$-algebra. Such an interior $P$-algebra $A_\gamma$ is called a source algebra of the block algebra $\O Gb$. By \cite[Corollary 3.5]{P0}, $A_\gamma$ and $\O Gb$ are Morita equivalent. In particular, given a simple $\K G$-module $M$ in the block $b$, $iM$ is a simple $\K\otimes_\O A_\gamma$-module, and the correspondence $$M\mapsto iM\leqno 3.5.1$$ gives a bijection between the set of isomorphism classes of simple $\K G$-modules in the block $b$ and the set of isomorphism classes of simple $\K\otimes_\O A_\gamma$-modules.

\bigskip We say that the block algebras $\O Gb$ and $\O G'b'$ are {\it basically Morita equivalent} if there is an $\O(G\times G')$-module with endopermutation source, inducing a Morita equivalence between the block algebras $\O Gb$ and $\O G'b'$. In this case, such a Morita equivalence is called a {\it basic Morita equivalence} between $\O Gb$ and $\O G'b'$ (see \cite{P3}).

\bigskip\noindent{\bf Proposition 3.6.}\quad {\it Let $G$ and $G'$ be finite groups and $b$ and $b'$ blocks of $G$ and $G'$.
Assume that the block algebras $\O Gb$ and $\O G'b'$ are basically Morita equivalent and that defect groups of $b$ and $b'$ are abelian. Then there is a Morita equivalence between the block algebras $\O Gb$ and $\O G'b'$ inducing a perfect isometry $I_\nu: \mathbb Z{\rm Irr}(G', b')\rightarrow \mathbb Z{\rm Irr}(G, b)$ such that $I_\nu^\K$ is compatible with a local system $\{I_T\}_{\{T ({\rm cyclic})\subset P\}}$.}

\bigskip\noindent{\it Proof.}\quad Since $\O Gb$ and $\O G'b'$ are basically Morita equivalent, by \cite[Corollary 7.4]{P3} they have isomorphic defect groups. For convenience, we assume without loss of generality that $b$ and $b'$ have a common defect group $P$. By \cite[6.9.3 and Corollary 7.4]{P3}, there are maximal local pointed groups $P_\gamma$ on $\O Gb$ and $P_{\gamma'}$ on $\O G'b'$ and an endopermutation $\O P$-module $N$ with vertex $P$, such that setting $(\O G)_\gamma=i(\O G)i$ and $(\O G')_{\gamma'}=i'(\O G')i'$ for some $i\in \gamma$ and some $i'\in \gamma'$, we have an interior $P$-algebra embedding $$f: (\O G)_\gamma\rightarrow {\rm End}_\O(N)\otimes_\O (\O G')_{\gamma'}.$$ Explicitly, $f$ is an injective $\O$-algebra homomorphism, the image of $f$ is equal to $f(i)\Big({\rm End}_\O(N)\otimes_\O (\O G')_{\gamma'}\Big)f(i)$, and $f$ preserves interior $P$-algebra structures (which is an interior $P$-algebra homomorphism in the sense of  \cite[Definition 3.1]{P0}); the interior $P$-algebra structure on ${\rm End}_\O(N)\otimes_\O (\O G')_{\gamma'}$ is determined diagonally by the interior $P$-algebra structures on ${\rm End}_\O(N)$ and $(\O G')_{\gamma'}$. We identify $(\O G)_\gamma$ as its image through $f$, so that $$(\O G)_\gamma=i\Big({\rm End}_\O(N)\otimes_\O (\O G')_{\gamma'}\Big)i.\leqno 3.6.1$$

Since $\O Gb$ and $\O G'b'$ are basically Morita equivalent, $(\O G)_\gamma$ and $(\O G')_{\gamma'}$ have the same numbers of isomorphism classes of simple $(\O G)_\gamma$- and $(\O G')_{\gamma'}$-modules. Clearly ${\rm End}_\O(N)\otimes_\O (\O G')_{\gamma'}$ and $(\O G')_{\gamma'}$ have the same numbers of isomorphism classes of simple $({\rm End}_\O(N)\otimes_\O (\O G')_{\gamma'})$- and $(\O G')_{\gamma'}$-modules. So $(\O G)_\gamma$ and ${\rm End}_\O(N)\otimes_\O (\O G')_{\gamma'}$ have the same numbers of isomorphism classes of simple $(\O G)_\gamma$- and $({\rm End}_\O(N)\otimes_\O (\O G')_{\gamma'})$-modules. By \cite[Theorem 9.9]{T}, the correspondence $$W\rightarrow i(N\otimes_\O W)\leqno 3.6.2$$ gives a Morita equivalence between the categories of $(\O G')_{\gamma'}$- and $(\O G)_\gamma$-modules. By composing the Morita equivalences between $\O G'b'$ and $(\O G')_{\gamma'}$ (see 3.5), between $(\O G')_{\gamma'}$ and $(\O G)_\gamma$ and between $(\O G)_\gamma$ and $\O Gb$ (see 3.5), we get a Morita equivalence $\Phi$ between $\O G'b'$ and $\O Gb$.

Note that $\Phi$ may be different from the Morita equivalence $\Psi$ between $\O G'b'$ and $\O Gb$ in the assumption. For a $\K$-algebra $A$, we denote by ${\cal M}(A)$ the set of isomorphism classes of simple $A$-modules.
The purpose of constructing $\Phi$ is to get the following commutative diagram $$\begin{array}{ccc}
  {\cal M}(\K G'b') & \rightarrow  & {\cal M}(\K Gb)  \\
  \downarrow &  & \downarrow \\
  {\cal M}(\K\otimes_\O(\O G')_{\gamma'}) & \rightarrow  & {\cal M}(\K\otimes_\O(\O G)_\gamma),
\end{array}$$ where the vertical arrows are obtained by the correspondence 3.5.1 applied to $\O Gb$ and $\O G'b'$, the top arrow $\phi$ is induced by $\Phi$ and the bottom arrow is determined by the correspondence 3.6.2. The diagram is a consequence of the Morita equivalence $\Phi$ in the last paragraph. The bijection $ {\cal M}(\K G'b') \rightarrow {\cal M}(\K Gb)  $ induced by $\Psi$ may not make the above diagram commutative.

We say that a local pointed group $Q_\delta$ on $\O G$ is associated to a Brauer pair $(Q, g)$, if $${\rm Br}_Q(f){\rm Br}_Q(\delta)={\rm Br}_Q(\delta)$$ and that a local pointed element $u_\alpha$ is associated to a Brauer element $(u, g)$ if the local pointed group $\langle u\rangle_\delta$ is associated to the Brauer pair $(\langle u\rangle, g)$. Let $(P, f)$ and $(P, f')$ be Brauer pairs, to which $P_\gamma$ and $P_{\gamma'}$ are associated. Let $M$ be a simple $\K Gb$-module and $\chi$ its ordinary character. Let $(u, h)$ be a Brauer element contained in $(P, f)$. Since $P$ is abelian, any local pointed element $u_\alpha$ associated to the Brauer element $(u, h)$ is contained in $P_\gamma$ and conversely any local pointed element $u_\alpha$ contained in $P_\gamma$ is associated to $(u, h)$. By 3.4.1, we have $$d_G^{(u,\, h)}(\chi)=\sum_{u_\alpha\in {\cal P}(u,\, h)}\chi^\alpha(u)\varphi_\alpha,\leqno 3.6.3$$
where ${\cal P}(u,\, h)$ denotes the set of all local pointed elements associated to the Brauer element $(u, h)$.
Let $M'$ be a simple $\K G'b'$-module corresponding to $M$ under the bijection $\phi$, and $\chi'$ its ordinary character. Let $(u, h')$ be the Brauer pair contained in $(P, f')$. Similarly, we have $$d_{G'}^{(u,\, h')}(\chi')=\sum_{u_{\alpha'}\in {\cal P}(u,\, h')}(\chi')^{\alpha'}(u)\varphi_{\alpha'}.\leqno 3.6.4$$

Given $u_\alpha\in {\cal P}(u,\, h)$, by \cite[Proposition 15.1]{T}, the intersection $\alpha\cap (\O G)_\gamma$ is still a local point of $\langle u\rangle$ on $(\O G)_\gamma$ and the correspondence $\alpha\mapsto \alpha\cap (\O G)_\gamma$ gives a bijection between ${\cal P}(u,\, h)$ and the set of local points of $\langle u\rangle$ on $(\O G)_\gamma$. Through the bijection, we identify ${\cal P}(u,\, h)$ as a set of local pointed elements on $(\O G)_\gamma$.
By \cite[7.6.2]{P3} applied to the equality 3.6.1, we get a bijection $${\cal P}(u,\, h)\rightarrow {\cal P}(u,\, h'),\, u_\alpha\mapsto u_{\alpha'}$$ where $u_\alpha$ and $u_{\alpha'}$ correspond to each other if and only if there are $j\in \alpha$ and $j'\in \alpha'$ such that $$j(\ell\otimes j')=j=(\ell\otimes j')j,$$ where $\ell$ is a primitive idempotent in the unique local point of $\langle u\rangle$ on ${\rm End}_\O (N)$. We are going to prove that $$\chi^\alpha(u)=\omega(u)(\chi')^{\alpha'}(u)\leqno 3.6.5$$ under the bijection, where $\omega$ denotes the character of the $\O \langle u\rangle$-module $\ell(N)$.

Since the $\K G$-module $M$ and the $\K G'$-module $M'$ correspond to each other under the bijection $\phi$, by the construction of $\Phi$ and 3.6.2, we have $$iM\cong i(N\otimes_\O i' M)$$ as $\K\otimes_\O (\O G)_\gamma$-modules. On the other hand, since $u_\alpha$ is contained in $P_\gamma$, there is $j''\in \alpha\cap (\O G)_\gamma$ such that $ij''=j''=j''i$. So $$j''M\cong j''(N\otimes_\O i' M)$$ as $(\K\otimes_\O (\O G)_\alpha)$-modules, where $(\O G)_\alpha=j''(\O G)j''$. We note that both $j$ and $j''$ are inside $\alpha\cap (\O G)_\gamma$ and thus that the idempotents $j''$ and $j$ are conjugate in $({\rm End}_\O(N)\otimes (\O G')_{\gamma'})^{\langle u \rangle}$. In particular, we have $$j''(N\otimes_\O i' M)\cong j(N\otimes_\O i' M)$$ as $\O  \langle u \rangle$-modules. Now, summarizing the above three isomorphisms, we conclude that the value of $\chi^\alpha$ at $u$ is equal to the value at $u$ of the character of the $\O  \langle u \rangle$-module $j(N\otimes_\O i' M)$.

Since $j(\ell\otimes j')=j=(\ell\otimes j')j$, we have $j(N\otimes_\O i' M)=j(\ell(N)\otimes_\O j' M)$. By \cite[5.6.3]{P4}, ${\rm Br}_{\langle u\rangle}^{{\rm End}_\O(N)\otimes_\O (\O G')_{\gamma'}}(\ell\otimes j'-j)=0$. So
the $\O \langle u \rangle$-module $(\ell\otimes j'-j)(N\otimes_\O i' M)$ is the direct sum of indecomposable $\O \langle u \rangle$-modules with vertex properly contained in $\langle u \rangle$, its character has the value 0 at $u$, and the value at $u$ of the character of the $\O  \langle u \rangle$-module $j(N\otimes_\O i' M)$ is equal to the value at $u$ of the character of the $\O  \langle u \rangle$-module $\ell(N)\otimes_\O i' M$, which is equal to $\omega(u)(\chi')^{\alpha'}(u)$. Until now, the equality 3.6.5 is proved.

The bijection ${\rm Irr}(G', b')\rightarrow {\rm Irr}(G, b), \,\chi'\mapsto \chi$ induced by $\Phi$ can be extended to a perfect isometry $I_\nu: \mathbb Z{\rm Irr}(G', b')\rightarrow \mathbb Z{\rm Irr}(G, b)$ by the $\mathbb Z$-linearity. We define a $\K$-linear isomorphism $$I_u: {\rm CF}(C_{G'}(u),\, h',\, \K)\rightarrow {\rm CF}(C_G(u),\, h,\, \K),$$ which sends $\varphi_{\alpha'}$ onto $\omega(u)\varphi_\alpha$ if local pointed elements $u_\alpha\in {\cal P}(u,\, h)$ and $u_{\alpha'}\in {\cal P}(u,\, h')$ correspond under the bijection ${\cal P}(u,\, h)\rightarrow {\cal P}(u,\, h')$. By the equalities 3.6.3, 3.6.4 and 3.6.5 we have $d_G^{(u,\,h)}\circ I_\nu^\K=I_u\circ d_{G'}^{(u,\,h')}$. In particular, $I_\nu^\K$ is compatible with a local system $\{I_u|u\in P\}$.

\bigskip\bigskip\noindent{\bf 4. Reduction}

\bigskip In this section, we give several lemmas for the reduction of the proof of {\bf Theorem}.

\bigskip\noindent{\bf 4.1.}\quad Let $G$ be a finite group and $b$ a block of $G$ with maximal $b$-Brauer pair $(P, f)$. Obviously $N_G(P, f)$ acts on the simple algebra $k\otimes_{\O Z(P)} \O C_G(P)f$ by the conjugation, where $Z(P)$ is the center of $P$. We denote by $\hat N_G(P, f)$ the set of all pairs $(x, a_x)$, where $x\in N_G(P, f)$ and $a_x$ is an invertible element of $k\otimes_{\O Z(P)} \O C_G(P)f$ such that $d^x=d^{a_x}$ for all $d\in k\otimes_{\O Z(P)} \O C_G(P)f$.
Obviously $\hat N_G(P, f)$ is a subgroup of the direct product $\hat N_G(P, f)\times (k\otimes_{\O Z(P)} \O C_G(P)f)^*$. The maps $k^*\rightarrow \hat N_G(P, f), \lambda\mapsto (1, \lambda)$ and $C_G(P)\rightarrow \hat N_G(P, f), x\mapsto (x, x(1\otimes f))$ are injective group homomorphisms. The image of $k^*$ is central in $\hat N_G(P, f)$, the image of $C_G(P)$ is normal in $\hat N_G(P, f)$ and they intersect trivially. We identify $k^*$ and $C_G(P)$ with their images through the two homomorphisms. The quotient of $\hat N_G(P, f)$ by $k^*$ is isomorphic to $N_G(P, f)$ and $\hat N_G(P, f)$ is a central extension of $N_G(P, f)$ by $k^*$. Set $E_G(P, f)=N_G(P, f)/C_G(P)$ and $\hat E_G(P, f)=\hat N_G(P, f)/C_G(P)$. The quotient group $E_G(P, f)$ is the inertial quotient of the block $b$ if $P$ is abelian. The homomorphism $k^*\rightarrow \hat N_G(P, f), \lambda\mapsto (1, \lambda)$ induces a new injective group homomorphism $k^*\rightarrow \hat E_G(P, f)$, and the quotient of $\hat E_G(P,\, f)$ by the image of $k^*$ is isomorphic to $E_G(P,\,f)$. So $\hat E_G(P, f)$ is a central extension of $E_G(P, f)$ by $k^*$.

\medskip\noindent{\bf 4.2.}\quad In this section, we always assume that $P$ is abelian and that the central extension $\hat E_G(P, f)$ splits.
The $N_G(P, f)$-conjugation induces actions of $E_G(P, f)$ on $P$. Set $K=P\rtimes E_G(P, f)$. We say that a perfect isometry $I_\nu: \mathbb Z{\rm Irr}(K)\rightarrow \mathbb Z {\rm Irr}(G, b)$ is {\it compatible with the $\ast$-structure} if $I_\nu(\lambda\ast\chi)=\lambda \ast I_\nu(\chi)$ for any $\chi\in {\rm Irr}(K)$ and any $N_G(P,\,f)$-stable character $\lambda$ of $P$, and that a Morita equivalence between $\O Gb$ and $\O K$ is a Morita equivalence compatible with the $\ast$-structure if the Morita equivalence induces a perfect isometry compatible with the $\ast$-structure.

\medskip\noindent{\bf 4.3.}\quad Let $N$ be a normal subgroup of $G$ and $h$ a block of $N$ covered by the block $b$. Let $H$ be the stabilizer of $h$ under the $G$-conjugation action. Then there is a suitable block $e$ of $H$ such that $b=\sum_{x}xex^{-1}$, where $x$ runs over a set of representatives of left cosets of $H$ in $G$, and such that $ee^x=0$ for any $x\in G-H$. The blocks $b$ and $e$ have common defect groups and the block algebra $\O Gb$ is isomorphic to the induction ${\rm Ind}^G_H(\O He)$ of the block algebra $\O He$ (see \cite[\S 16]{T}). In particular, there is a Morita equivalence between the block algebras $\O Gb$ and $\O He$ induced by an indecomposable $\O(G\times H)$-module with trivial sources. We adjust the choice of $P$, so that $P$ is also a defect group of the block $e$. By \cite[7.6.5]{P3}, $\hat E_G(P, f)$ is isomorphic to $\hat E_H(P, f')$ as central extensions for some maximal $e$-Brauer pair $(P, f')$. In particular, $E_G(P, f)$ is isomorphic to $E_H(P, f')$. We identify $E_G(P, f)$ and $E_H(P, f')$.
Since the central extension $\hat E_G(P, f)$ is assumed to split, $\hat E_H(P, f')$ splits.

\bigskip\noindent{\bf Lemma 4.4.}\quad {\it Keep the notation and the assumptions of 4.2 and 4.3. Assume that there is a Morita equivalence between the algebras $\O He$ and $\O K$ compatible with the $\ast$-structure. There is a Morita equivalence between the block algebras $\O Gb$ and $\O K$ compatible with the $\ast$-structure. }

\medskip\noindent{\it Proof.}\quad Since there is a Morita equivalence between the block algebras $\O Gb$ and $\O He$ induced by an indecomposable $\O(G\times H)$-module with trivial sources, by Propositions 2.7 and 3.6, there is a Morita equivalence inducing a perfect isometry $I_\nu: \mathbb Z{\rm Irr}(H, e)\rightarrow \mathbb Z {\rm Irr}(G, b)$ such that $I_\nu(\lambda\ast \chi)=\lambda\ast I_\nu(\chi)$ for any $\chi\in {\rm Irr}(H, e)$ and any $N_G(P, f)$-stable character of $P$.
The composition of $I_\nu$ and the perfect isometry induced by the Morita equivalence between $\O He$ and $\O K$ is a perfect isometry compatible with the $\ast$-structure, which is induced by the composition of the latter Morita equivalence between $\O Gb$ and $\O He$ and the Morita equivalence between $\O He$ and $\O K$.

\bigskip\noindent{\bf Lemma 4.5.}\quad {\it Keep the notation and the assumption of 4.2 and 4.3. Set $Q=P\cap N$. Then $Q$ is a defect group of the block $h$. Assume that $H=G$ and that the block $h$ is nilpotent. Then there is a finite group $L$ such that

\smallskip\noindent{\bf 4.5.1.}\quad $Q$ is normal in $L$ and $P$ is a Sylow $p$-subgroup of $L$,

\smallskip\noindent{\bf 4.5.2.}\quad the quotient groups $G/N$ and $L/Q$ are isomorphic, and

\smallskip\noindent{\bf 4.5.3.}\quad there are a $p'$-central extension $\tilde L$ of $L$ and a block $\ell$ of $\tilde L$ with defect group $P$, such that the block algebras $\O Gb$ and $\O \tilde L\ell$ are basically Morita equivalent.
}

\medskip\noindent{\it Proof.} This lemma follows from \cite[Theorem 1.12]{KP}.

\bigskip By \cite[7.6.5]{P3}, $\hat E_G(P, f)$ is isomorphic to $\hat E_{\tilde L}(P, f')$ as central extensions for some maximal $\ell$-Brauer pair $(P, f')$. In particular, $E_G(P, f)$ is isomorphic to $E_{\tilde L}(P, f')$. We identify $E_G(P, f)$ and $E_{\tilde L}(P, f')$.
Since the central extension $\hat E_G(P, f)$ is assumed to split, $\hat E_{\tilde L}(P, f')$ splits.

\bigskip\noindent{\bf Lemma 4.6.}\quad {\it Keep the notation and the assumptions in Lemma 4.5. Assume that there is a Morita equivalence between the algebras $\O {\tilde L}\ell$ and $\O K$ compatible with the $\ast$-structure. Then there is a Morita equivalence between the algebras $\O G b$ and $\O K$ compatible with the $\ast$-structure. }

\medskip\noindent{\it Proof.} One can prove the lemma, using the same proof method for Lemma 4.4.

\bigskip\bigskip\noindent{\bf 5. Extensions of Morita equivalences}

\bigskip In this section, we mainly give an extension of Morita equivalences compatible with the $\ast$-structure (see Proposition 5.9 below).

\bigskip\noindent{\bf 5.1.}\quad Let $G$ be a finite group, $H$ a normal subgroup of $G$ with $p$-power index, and $b$ a $G$-stable block of $H$. Then $b$ is also a block of $G$. We always assume that $b$ as a block of $G$ has abelian defect group $P$. Set $Q=P\cap H$. By \cite[Proposition 5.3]{KP}, $Q$ is a defect group of the block $b$ of $H$ and $G$ is equal to $PH$.
Let $(P, f)$ be a maximal Brauer pair of the block $b$ of $G$.
Denote by $T$ the subgroup of all $E_G(P, f)$-fixed elements in $P$ and by $R$ the commutator $[P, N_G(P, f)]$. We have $P=T\times R$, $R=[P, N_H(P, f)]\subset H$, and $R$ is contained in $Q$. In particular, we have $G=TH$. Set $E=E_G(P, f)$, $K=P\rtimes E$ and $L=Q\rtimes E$. Then $K=LT$.

\bigskip\noindent{\bf Lemma 5.2.}\quad {\it Keep the notation and the assumption of 5.1. Any simple module in the block $b$ of $H$ is extendible to $G$ and any irreducible ordinary character in the block $b$ of $H$ is extendible to $G$.}

\medskip\noindent{\it Proof.}\quad See \cite[Theorem]{Kr}.

\bigskip  We label characters in ${\rm Irr}(K)$ as $\zeta_j$, $j\in J$, and choose a subset $I\subset J$, such that all these restrictions $(\zeta_i)_L, \,i\in I$, are exactly all characters in ${\rm Irr}(L)$. For any $i\in I$, set $\phi_i=(\zeta_i)_L$. Then ${\rm Irr}(L)=\{\phi_i| i\in I\}$.

\bigskip\noindent{\bf Lemma 5.3.}\quad {\it Keep the notation and the assumption of 5.1. Assume that $I_\nu: \mathbb Z{\rm Irr}(K)\rightarrow \mathbb Z{\rm Irr}(G, b)$ is a perfect isometry compatible with the $\ast$-structure.
Then there is a perfect isometry $I_{\nu'}: \mathbb Z{\rm Irr}(L)\rightarrow \mathbb Z{\rm Irr}(H, b)$
compatible the $\ast$-structure such that the diagram $$\begin{array}{ccc}
  \mathbb Z{\rm Irr}(K)  & \xrightarrow[I_\nu]{}  & \mathbb Z{\rm Irr}(G, b)  \\
  \downarrow &  & \downarrow \\
  \mathbb Z{\rm Irr}(L) & \xrightarrow{I_{\nu'}}  & \mathbb Z{\rm Irr}(H, b)
\end{array}\leqno 5.3.1$$ is commutative, where the two down arrows are determined by the restriction of characters. Moreover, $\nu'$ has a suitable extension $\hat\nu'$ to $(H\times L)\Delta(P)$ such that $\nu={\rm Ind}^{G\times K}_{(H\times L)\Delta(P)}(\hat\nu')$.

}

\medskip\noindent{\it Proof.}\quad
Since $I_\nu: \mathbb Z{\rm Irr}(K)\rightarrow \mathbb Z{\rm Irr}(G, b)$ is a perfect isometry compatible with the $\ast$-structure, $$I_\nu(\lambda\ast\chi)=\lambda\ast I_\nu(\chi)\leqno 5.3.2$$ for any $\chi\in \mathbb Z{\rm Irr}(K)$ and any $N_G(P, f)$-stable character $\lambda$ of $P$.

For any $i, j\in I$, there exist $\epsilon_i, \epsilon_j\in \{\pm 1\}$ and  $\chi_i, \chi_j\in {\rm Irr}(G, b)$, such that $\epsilon_i\chi_i=I_\nu(\zeta_i)$ and $\epsilon_j\chi_j=I_\nu(\zeta_j)$. By Lemma 5.2, the restrictions $\psi_i$ and $\psi_j$ of $\chi_i$ and $\chi_j$ to $H$ are irreducible. Suppose $\psi_i=\psi_j$ for different $i$ and $j$. Set $S=T\cap H$. The quotient group $G/H$ is naturally isomorphic to $T/S$. Through the isomorphism, any linear character $\lambda$ of $T$ with kernel containing $S$ can be regarded as a character of $G$ with kernel containing $H$. Since $\psi_i=\psi_j$,
there exists a linear character $\lambda$ of $T$ with kernel containing $S$, such that $\chi_i=\chi_j\lambda$. Obviously $\lambda$ is $E_G(P, f)$-stable, and by \cite[Corollary]{BP} $\chi_j\lambda$ is equal to $\lambda \ast \chi_j$. By 5.3.2, $\lambda\epsilon_j\zeta_j=\epsilon_i\zeta_i$, $\epsilon_i=\epsilon_j$ and $\zeta_i=\zeta_j$. Thus we have $\phi_i=\phi_j$. That contradicts the characters $\phi_i$ and $\phi_j$ being different for different $i$ and $j$. On the other hand, denoting by ${\rm Irr}(T/S)$ the set of irreducible ordinary characters of $T$ with kernel containing $S$, by 5.3.2, we have $\mathrm{Irr}(G, b)=
\{\lambda\chi_i\mid i\in I; \lambda\in {\rm Irr}(T/S)\}$.
Therefore $\mathrm{Irr}(H, b)=
\{\psi_i|i\in I\}$.

By 5.3.2, the perfect isometry $I_\nu$ sends $\lambda\zeta_i$ onto $\epsilon_i\lambda\chi_i$ for any $i\in I$ and any $\lambda\in {\rm Irr}(T/S)$. So $\nu=\sum_{i\in I}\sum_{\lambda\in {\rm Irr}(T/S)}\epsilon_i(\lambda\chi_i)\times (\lambda\zeta_i)^\ast$. Set $\nu'=\sum_{i\in I}\epsilon_i(\chi_i\times \zeta_i^\ast)$. Then $\nu={\rm Ind}^{G\times K}_{(H\times L)\Delta(P)}(\nu')$.
Considering $\nu'$ as a generalized character on $H\times L$ by restriction, the map $I_{\nu'}: \mathbb Z{\rm Irr}(L)\rightarrow \mathbb Z {\rm Irr}(H, b)$ sends $\phi_i$ onto $\epsilon_i\psi_i$ for any $i\in I$ and it is a bijective isometry. Moreover, it is trivial to see that $I_\nu$ and $I_{\nu'}$ satisfies the diagram 5.3.1. By the following lemma, the proof is done.

\bigskip\noindent{\bf Lemma 5.4.}\quad {\it Keep the notation and the assumptions in Lemma 5.3 and its proof. Then the map $I_{\nu'}$ is a perfect isometry compatible with the $\ast$-structure.
}

\medskip\noindent{\it Proof.}\quad Let $(x, f_{\langle x\rangle})$ be a $b$-Brauer element contained in $(P, f)$. The pair $(P, f)$ is also a maximal $b_{\langle x\rangle}$-Brauer pair and the quotient group $N_{C_G(x)}(P, f)/C_G(P)$ is isomorphic to $C_E({\langle x\rangle})$. Set $K_{\langle x\rangle}=P\rtimes C_E({\langle x\rangle})$. By Proposition 2.7, there is a $\K$-linear isomorphism $I_{\langle x\rangle}: {\rm CF}_{p'}(K_{\langle x\rangle}, \K)\rightarrow {\rm CF}_{p'}(C_G(x), b_{\langle x\rangle}, \K)$,
such that $$I_{\langle x\rangle}\circ d_K^x=d_G^{(x,\, b_{\langle x\rangle})}\circ I_\nu^\K.\leqno 5.4.1$$

Denote by ${\rm Res}^G_H$ the restriction map from ${\rm CF}(G, b, \K)\rightarrow {\rm CF}(H, b, \K)$ and by ${\rm Res}^G_{H,\,p'}$ the restriction map from ${\rm CF}_{p'}(G, b, \K)\rightarrow {\rm CF}_{p'}(H, b, \K)$. By Lemma 5.2 ${\rm Res}^G_H$ is surjective and ${\rm Res}^G_{H,\,p'}$ is a $\K$-linear isomorphism.

Assume that $x$ lies in $Q$. Note that $b_{\langle x\rangle}$ is a block of $C_H(x)$ with defect group $Q$. We set $$L_{\langle x\rangle}=Q\rtimes C_E({\langle x\rangle})\,\,\mbox{\rm and}\,\, I_{\langle x\rangle}'={\rm Res}^{C_G(x)}_{C_H(x),\,p'}\circ I_{\langle x\rangle}\circ ({\rm Res}^{K_{\langle x\rangle}}_{L_{\langle x\rangle},\,p'})^{-1}.\leqno 5.4.2$$ We claim that $d_H^{(x,\, b_{\langle x\rangle})}\circ I_{\nu'}^\K=I_{\langle x\rangle}'\circ d_L^x$. Indeed, it is easy to prove that the equalities $${\rm Res}^{C_G(x)}_{C_H(x),\,p'}\circ d_G^{(x,\, b_{\langle x\rangle})}=d_H^{(x,\, b_{\langle x\rangle})}\circ {\rm Res}^G_H\,\,\mbox{\rm and}\,\,{\rm Res}^{K_{\langle x\rangle}}_{L_{\langle x\rangle}}\circ d_K^x=d_L^x\circ {\rm Res}^K_L\leqno 5.4.3$$ hold. Since
$$\begin{array}{ccl}
  d_H^{(x,\, b_{\langle x\rangle})}\circ I_{\nu'}^\K\circ {\rm Res}^K_L & = & d_H^{(x,\, b_{\langle x\rangle})}\circ {\rm Res}^G_H\circ I_\nu^\K \\
   & = & {\rm Res}^{C_G(x)}_{C_H(x),\,p'}\circ d_G^{(x,\, b_{\langle x\rangle})}\circ I_\nu^\K\\
   &=& {\rm Res}^{C_G(x)}_{C_H(x),\,p'}\circ I_{\langle x\rangle}\circ d_K^x\\
   &=& I_{\langle x\rangle}'\circ {\rm Res}^{K_{\langle x\rangle}}_{L_{\langle x\rangle},\,p'} \circ d_K^x\\
   &=& I_{\langle x\rangle}'\circ d_L^x\circ {\rm Res}^K_L
\end{array}$$ where the first equality is obtained by 5.3.1, the second by 5.4.3, the third by 5.4.1, the fourth by 5.4.2 and the fifth by 5.4.3,
the claim is done. Therefore $I_{\nu'}^\K$ is compatible with a local system $\{I_{\langle x\rangle}'|\mbox{$x\in P$}\}$ in the sense of \cite[Definition 4.3]{B1}. By Proposition 2.7, $I_{\nu'}^\K$ is compatible with the $\ast$-structure.

In order to prove that $I_{\nu'}$ is a perfect isometry, by \cite[Lemma 4.5]{B1}, it remains to prove that $I_{\langle x\rangle}'$ maps ${\rm CF}_{p'}(L_{\langle x\rangle}, \O)$ onto ${\rm CF}_{p'}(C_H(x), b_{\langle x\rangle}, \O)$. Since $I_{\langle x\rangle}'={\rm Res}^{C_G(x)}_{C_H(x),\,p'}\circ I_{\langle x\rangle}\circ ({\rm Res}^{K_{\langle x\rangle}}_{L_{\langle x\rangle},\,p'})^{-1}$, it suffices to show that $I_{\langle x\rangle}$ maps ${\rm CF}_{p'}(K_{\langle x\rangle}, \O)$ onto ${\rm CF}_{p'}(C_G(x), b_{\langle x\rangle}, \O)$.
This follows from Lemma 2.8.

\bigskip Set $L'=R\rtimes E$. Then $L$ is equal to the direct product $L'\times T'$, where $T'$ is the subgroup of all $E$-fixed elements in $Q$. Let $I_{\psi'}: \mathbb Z {\rm Irr}(L')\rightarrow \mathbb Z {\rm Irr}(L')$ be a perfect isometry compatible with the $\ast$-structure. Since any irreducible ordinary character of $L$ is of the form $\phi\times \lambda$, where $\phi\in {\rm Irr}(L')$ and $\lambda\in {\rm Irr}(T')$, we extend $I_{\psi'}$ to a $\mathbb Z$-linear map $I_{\psi}: \mathbb Z {\rm Irr}(L)\rightarrow \mathbb Z {\rm Irr}(L)$ sending $\phi\times \lambda$ onto $I_{\psi'}(\phi)\times \lambda$ for any $\phi\in {\rm Irr}(L')$ and any $\lambda\in {\rm Irr}(T')$. Set $\tilde\psi'=\psi'\times 1_{\Delta(T')}$, since $(L'\times L')\Delta(T')=(L'\times L')\times \Delta(T')$.

\bigskip\noindent{\bf Lemma 5.5.}\quad {\it Keep the notation and the assumption as above. Then $I_{\psi}$ is a perfect isometry compatible with the $\ast$-structure, $\psi={\rm Ind}^{L\times L}_{(L'\times L')\Delta(T')}(\tilde\psi')$, and the diagram $$\begin{array}{ccc}
  \mathbb Z{\rm Irr}(L) & \xrightarrow[I_{\psi}]{}  & \mathbb Z{\rm Irr}(L)  \\
  \downarrow &  & \downarrow \\
  \mathbb Z{\rm Irr}(L') & \xrightarrow{I_{\psi'}}  & \mathbb Z{\rm Irr}(L')
\end{array}$$ is commutative, where the downarrows are induced by the restriction of characters. }

\medskip\noindent{\it Proof.}\quad The last statement of the lemma is obvious. Since
$$\psi=\sum_{\phi\in {\rm Irr}(L')}\sum_{\lambda\in {\rm Irr}(T')} (I_{\psi'}(\phi)\times \lambda)\times (\phi^*\times \lambda^*)=\psi'\times \Big(\sum_{\lambda\in {\rm Irr}(T')}\lambda\times\lambda^\ast\Big),$$ it is easy to see $\psi={\rm Ind}^{L\times L}_{(L'\times L')\Delta(T')}(\tilde\psi')$. It remains to show that $I_{\psi}$ is a perfect isometry compatible with the $\ast$-structure.

Since any irreducible ordinary character of $L$ is of the form $\phi\times \lambda$, where $\phi\in {\rm Irr}(L')$ and $\lambda\in {\rm Irr}(T')$, there is an obvious group isomorphism $$\mathbb Z {\rm Irr}(L)\cong \mathbb Z {\rm Irr}(L')\bigotimes _{\mathbb Z} \mathbb Z {\rm Irr}(T').$$ We identify $\mathbb Z {\rm Irr}(L)$ and $\mathbb Z {\rm Irr}(L')\bigotimes _{\mathbb Z} \mathbb Z {\rm Irr}(T')$ through the isomorphism. Since $Q=R\times T'$, we also have a $\K$-linear isomorphism $${\rm CF}_{p'}(Q,\, \K)\cong {\rm CF}_{p'}(R,\, \K)\bigotimes_\K {\rm CF}_{p'}(T',\, \K), 1_Q^\circ\mapsto 1_R^\circ\otimes 1_{T'}^\circ$$ where $1_Q^\circ$ denotes the unique irreducible Brauer character on $Q$. We identify ${\rm CF}_{p'}(Q,\, \K)$ and ${\rm CF}_{p'}(R,\, \K)\bigotimes_\K {\rm CF}_{p'}(T',\, \K)$. With the identifications as above, we have $$I_{\psi}=I_{\psi'}\otimes {\rm Id}_{\mathbb Z {\rm Irr}(T')}\leqno 5.5.1$$ where ${\rm Id}_{\mathbb Z {\rm Irr}(T')}$ is the identity map on $\mathbb Z {\rm Irr}(T')$, and $$d^{xu}_L=d_{L'}^x\otimes d_{T'}^u\leqno 5.5.2$$ for any $x\in R$ and any $u\in T'$. On the other hand, since $I_{\psi'}$ is a perfect isometry compatible with the $\ast$-structure, by Proposition 2.7, for any $x\in R$, there is a $\K$-linear map $I_x: {\rm CF}_{p'}(C_{L'}(x),\, \K)\rightarrow {\rm CF}_{p'}(C_{L'}(x),\, \K)$ such that $$d_{L'}^x\circ I_{\psi'}^\K=d_{L'}^x\circ I_x. \leqno 5.5.3$$
Now set $I_{xu}=I_x\otimes_\K {\rm Id}_{{\rm CF}_{p'}(T',\, \K)}$, where ${\rm Id}_{{\rm CF}_{p'}(T',\, \K)}$ is the identity map on ${\rm CF}_{p'}(T',\, \K)$. By 5.5.1, 5.5.2 and 5.5.3, we have $I_{xu}\circ d^{xu}_L= d^{xu}_L\circ I_{\psi}^\K$. By Proposition 2.7, $I_\psi$ is compatible with the $\ast$-structure.
Since $I_{\psi'}$ is a perfect isometry, by Lemma 2.8, $I_x$ maps ${\rm CF}_{p'}(C_{L'}(x),\, \O)$ onto ${\rm CF}_{p'}(C_{L'}(x),\, \O)$. Thus for any $x\in R$ and $u\in T'$, $I_{xu}$ maps ${\rm CF}_{p'}(C_{L}(xu),\, \O)$ onto ${\rm CF}_{p'}(C_{L}(xu),\, \O)$. By \cite[Lemma 4.5]{B1}, $I_{\psi}$ is a perfect isometry.

\bigskip\noindent{\bf Lemma 5.6.}\quad {\it Keep the notation and the assumption of 5.1. Assume that $E$ is cyclic and acts freely on $R-\{1\}$ and that the map $I_\psi: \mathbb Z {\rm Irr}(L)\rightarrow \mathbb Z {\rm Irr}(L)$ is a perfect isometry compatible with the $\ast$-structure. Then $\psi$ has a suitable extension $\tilde\psi$ to $(L\times L)\Delta(K)$, such that setting $\hat\psi={\rm Ind}^{K\times K}_{(L\times L)\Delta(K)}(\tilde\psi)$, the map $I_{\hat\psi}: \mathbb Z {\rm Irr}(K)\rightarrow \mathbb Z {\rm Irr}(K)$ is a perfect isometry compatible with the $\ast$-structure and makes the following diagram commutative $$\begin{array}{ccc}
  \mathbb Z{\rm Irr}(K) & \xrightarrow[I_{\hat\psi}]{}  & \mathbb Z{\rm Irr}(K)  \\
  \downarrow &  & \downarrow \\
  \mathbb Z{\rm Irr}(L) & \xrightarrow{I_\psi}  & \mathbb Z{\rm Irr}(L)
\end{array}\leqno 5.6.1$$ where the downarrows are induced by the restriction of characters.}

\medskip\noindent{\it Proof.}\quad By Lemma 5.3, there is a generalized character $\psi'$ of $(L'\times L')\Delta(T')$ such that
$\psi={\rm Ind}^{L\times L}_{(L'\times L')\Delta(T')}(\psi')$ and such that, restricting $\psi'$ to $L'\times L'$, the map $I_{\psi'}: \mathbb Z {\rm Irr}(L')\rightarrow \mathbb Z {\rm Irr}(L')$ is a perfect isometry compatible with the $\ast$-structure and makes the following diagram commutative $$\begin{array}{ccc}
  \mathbb Z{\rm Irr}(L) & \xrightarrow[I_\psi]{}  & \mathbb Z{\rm Irr}(L)  \\
  \downarrow &  & \downarrow \\
  \mathbb Z{\rm Irr}(L') & \xrightarrow{I_{\psi'}}  & \mathbb Z{\rm Irr}(L')
\end{array}$$
where the two down arrows are determined by the restriction of characters.

Set $\hat\psi'=\psi'\times 1_{\Delta(T')}$ and $\psi''={\rm Ind}^{L\times L}_{(L'\times L')\Delta(T')}(\hat\psi')$.
By Lemma 5.5 applied to $I_{\psi'}$ and the group $L$, the map $I_{\psi''}: \mathbb Z {\rm Irr}(L)\rightarrow \mathbb Z {\rm Irr}(L)$ is a perfect isometry compatible with the $\ast$-structure and the diagram $$\begin{array}{ccc}
  \mathbb Z{\rm Irr}(L) & \xrightarrow[I_{\psi''}]{}  & \mathbb Z{\rm Irr}(L)  \\
  \downarrow &  & \downarrow \\
  \mathbb Z{\rm Irr}(L') & \xrightarrow{I_{\psi'}}  & \mathbb Z{\rm Irr}(L')
\end{array}$$ is commutative  where the downarrows are induced by the restriction of characters. Note $(L'\times L')\Delta(K)=(L'\times L')\times \Delta(T)$.
Set $\tilde\psi'=\psi'\times 1_{\Delta(T)}$ and $\psi'''={\rm Ind}^{K\times K}_{(L'\times L')\Delta(T)}(\tilde \psi')$. Similarly,
the map $I_{\psi'''}: \mathbb Z {\rm Irr}(K)\rightarrow \mathbb Z {\rm Irr}(K)$ is a perfect isometry compatible with the $\ast$-structure and the diagram $$\begin{array}{ccc}
  \mathbb Z{\rm Irr}(K) & \xrightarrow[I_{\psi'''}]{}  & \mathbb Z{\rm Irr}(K)  \\
  \downarrow &  & \downarrow \\
  \mathbb Z{\rm Irr}(L') & \xrightarrow{I_{\psi'}}  & \mathbb Z{\rm Irr}(L')
\end{array}$$ is commutative, where the downarrows are induced by the restriction of characters.

Notice that the restriction of $\tilde\psi'$ to $(L'\times L')\Delta(T')$ is equal to $\hat\psi'$. Set $\hat \psi''={\rm Ind}_{(L'\times L')\Delta(T)}^{(L\times L)\Delta(T)}(\tilde \psi')$. The restriction to $L\times L$ of $\hat\psi''$ is equal to $\psi''$ and ${\rm Ind}^{K\times K}_{(L\times L)\Delta(T)}(\hat\psi'')$ is equal to $\psi'''$.
Moreover, $I_{\psi''}$ and $I_{\psi'''}$ are perfect isometries compatible with the $\ast$-structure and making the following diagram commutative $$\begin{array}{ccc}
  \mathbb Z{\rm Irr}(K) & \xrightarrow[I_{\psi'''}]{}  & \mathbb Z{\rm Irr}(K)  \\
  \downarrow &  & \downarrow \\
  \mathbb Z{\rm Irr}(L) & \xrightarrow{I_{\psi''}}  & \mathbb Z{\rm Irr}(L)
\end{array}$$ is commutative, where the downarrows are induced by the restriction of characters.
Therefore, to extend $I_\psi$ to a perfect isometry from $\mathbb Z {\rm Irr}(K)$ to itself compatible with the $\ast$-structure satisfying the diagram 5.6.1, is equivalent to extend the composition $I_\psi\circ I_{\psi''}^{-1}$ to a perfect isometry from $\mathbb Z {\rm Irr}(K)$ to itself compatible with the $\ast$-structure satisfying the diagram 5.6.1.

We replace $I_\psi$ by the composition $I_\psi\circ I_{\psi''}^{-1}$. Then, for any $\phi\in {\rm Irr}(L')$, there is some suitable $\lambda_\phi\in {\rm Irr}(T')$ such that $I_\psi(\phi\times 1_{T'})=\phi\times \lambda_\phi$. We claim that $\lambda_\phi=\lambda_{\phi'}$ for any $\phi, \phi'\in {\rm Irr}(L')$. In order to prove the claim, we take $\phi'$ to be $1_{L'}$. Since the degree of $\phi$ is coprime to $p$, by \cite[Theorem 5.16]{N}, there is some $u\in R-\{1\}$ such that $\phi(u)\neq 0$. Note that $C_{L}(uv)=Q$ since $E$ acts freely on $R-\{1\}$. For any $v\in T'$, since $I_\psi$ is a perfect isometry compatible with the $\ast$-structure, there is a $\K$-linear isomorphism $I_{\langle uv\rangle}: {\rm CF}_{p'}(Q, \K)\rightarrow  {\rm CF}_{p'}(Q, \K)$ such that $$I_{\langle uv\rangle}\circ d_L^{uv}=d_L^{uv}\circ I_\psi^\K.$$ Applying both sides of the equality to $\phi\times 1_{T'}$ and $\phi'\times 1_{T'}$, we get
$I_{\langle uv\rangle}(\phi(u)1_Q^\circ)=\phi(u)\lambda_{\phi}(v)1_Q^\circ$, $I_{\langle uv\rangle}(\phi'(u)1_Q^\circ)=\phi'(u)\lambda_{\phi'}(v)1_Q^\circ$, and thus $\lambda_\phi(v)=\lambda_{\phi'}(v)$ for any $v\in T'$. The claim is done.

We set $\pi=\lambda_{1_{L'}}$ and fix an extension $\varpi$ of $\pi$ to $T$. We extend $I_\psi$ to a $\mathbb Z$-linear map $I_{\hat\psi}: \mathbb Z {\rm Irr}(K)\rightarrow \mathbb Z {\rm Irr}(K)$ sending $\phi\times \lambda$ onto $\phi\times \varpi\lambda$ for any $\phi\in {\rm Irr}(L')$ and any $\lambda\in {\rm Irr}(T)$. It is trivial to verify that $I_{\hat\psi}$ is a perfect isometry compatible with the $\ast$-structure, that $\hat\psi={\rm Ind}^{K\times K}_{(L\times L)\Delta(T)}(\tilde\psi)$ for some suitable extension $\tilde\psi$ of $\psi$ to $(L\times L)\Delta(T)$, and that $I_\psi$ and $I_{\hat\psi}$ satisfy the diagram 5.6.1.

\bigskip\noindent{\bf Proposition 5.7.}\quad {\it Keep the notation and the assumption of 5.1 and assume that $E$ is cyclic and acts freely on $R-\{1\}$ and that $I_{\nu'}: \mathbb Z{\rm Irr}(L)\rightarrow \mathbb Z{\rm Irr}(H, b)$ is a perfect isometry compatible with the $\ast$-structure. Then there is a suitable extension $\tilde\nu'$ of $\nu'$ to $(H\times L)\Delta(P)$ such that setting $\nu={\rm Ind}^{G\times K}_{(H\times L)\Delta(P)}(\tilde\nu')$, the map $I_{\nu}: \mathbb Z{\rm Irr}(K)\rightarrow \mathbb Z {\rm Irr}(G, b)$ is a perfect isometry compatible with the $\ast$-structure and the diagram is commutative

$$\begin{array}{ccc}
  \mathbb Z{\rm Irr}(K) & \xrightarrow[I_\nu]{}  & \mathbb Z{\rm Irr}(G, b)  \\
  \downarrow &  & \downarrow \\
  \mathbb Z{\rm Irr}(L) & \xrightarrow{I_{\nu'}}  & \mathbb Z{\rm Irr}(H, b)
\end{array}$$
where the down arrows are induced by the restriction of characters.}

\medskip\noindent{\it Proof.} Denote by $\ell_G(b)$ the number of isomorphism classes of simple modules in the block $b$ of $G$. Since $\ell_H(b)$ is equal to $|E_G(P, f)|$ and $\ell_G(b)$ is equal to $\ell_H(b)$ (see Lemma 5.2), $\ell_G(b)$ is equal to $E_G(P, f)$.
By \cite[Theorem]{W2}, there is a perfect isometry $\mathbb Z{\rm Irr}(K)\rightarrow \mathbb Z{\rm Irr}(G, b)$ compatible with the $\ast$-structure. By Lemma 5.3, there is a perfect isometry $I_{\mu'}: \mathbb Z{\rm Irr}(L)\rightarrow \mathbb Z{\rm Irr}(H, b)$ compatible with the $\ast$-structure, $\mu={\rm Ind}^{G\times K}_{(H\times L)\Delta(P)}(\tilde\mu')$ for some extension $\tilde\mu'$ of $\mu'$ to $(H\times L)\Delta(P)$, and the diagram $$\begin{array}{ccc}
  \mathbb Z{\rm Irr}(K) & \xrightarrow[I_\mu]{}  & \mathbb Z{\rm Irr}(G, b)  \\
  \downarrow &  & \downarrow \\
 \mathbb Z{\rm Irr}(L)  & \xrightarrow{I_{\mu'}}  & \mathbb Z{\rm Irr}(H, b)
\end{array}$$ is commutative, where the two down arrows are determined by the restriction of characters.

Obviously $I_{\mu'}^{-1}\circ I_{\nu'}$ is a perfect isometry compatible with the $\ast$-structure. Suppose $I_{\eta'}=I_{\mu'}^{-1}\circ I_{\nu'}$. By Lemma 5.6, $\eta'$ has an extension $\tilde\eta'$ to $(L\times L)\Delta(K)$, such that setting $\eta={\rm Ind}^{K\times K}_{(L\times L)\Delta(K)}(\eta')$, the map $I_\eta: \mathbb Z {\rm Irr}(K)\rightarrow \mathbb Z {\rm Irr}(K)$ is a perfect isometry compatible with the $\ast$-structure and makes the following diagram commutative $$\begin{array}{ccc}
  \mathbb Z{\rm Irr}(K) & \xrightarrow[I_\eta]{}  & \mathbb Z{\rm Irr}(K)  \\
  \downarrow &  & \downarrow \\
  \mathbb Z{\rm Irr}(L) & \xrightarrow{I_{\eta'}}  & \mathbb Z{\rm Irr}(L)
\end{array}$$ where the downarrows are induced by the restriction of characters. Now we compose $I_\mu$ and $I_\eta$ and suppose $I_\nu=I_\mu\circ I_\eta$. Then $I_\nu$ is a perfect isometry compatible with the $\ast$-structure, and it makes commutative the diagram in the proposition. Notice that $I_\nu$ uniquely determines $I_{\nu'}$. By Lemma 5.3, $\nu'$ has an extension $\tilde\nu'$ to $(H\times L)\Delta(P)$ and $\nu={\rm Ind}^{G\times K}_{(H\times L)\Delta(P)}(\tilde\nu')$.

\medskip\noindent{\bf 5.8.}\quad
Assume that there are perfect isometries \begin{center}$I_\nu: \mathbb Z{\rm Irr}(K)\rightarrow \mathbb Z{\rm Irr}(G, b)$ and $I_{\nu'}: \mathbb Z{\rm Irr}(L) \rightarrow \mathbb Z{\rm Irr}(H, b)$,\end{center} such that $\nu={\rm Ind}^{G\times K}_{(H\times L)\Delta(P)}(\tilde\nu')$ for a suitable extension $\tilde\nu'$ of $\nu'$ to $(H\times L)\Delta(P)$, and such that $I_\nu$ and $I_{\nu'}$ make the following diagram commutative $$\begin{array}{ccc}
  \mathbb Z{\rm Irr}(K) & \rightarrow  & \mathbb Z{\rm Irr}(G, b)  \\
  \downarrow &  & \downarrow \\
  \mathbb Z{\rm Irr}(L) & \rightarrow  & \mathbb Z{\rm Irr}(H, b)
\end{array}\leqno 5.8.1$$ where the downarrows are induced by the restrictions of characters. By \cite[Theorem 1.5]{B1}, the perfect isometries $I_\nu$ and $I_{\nu'}$ induce algebra isomorphisms \begin{center} $\pi: Z(\O L)\cong Z(\O Hb)$ and $\tau: Z(\O K)\cong Z(\O Gb)$.   \end{center} Obviously $Z(\O Gb)$ and $Z(\O K)$ are $P/Q$-graded algebras with the $\bar u$-components $Z(\O Gb)\cap \O Hu$ and $Z(\O K)\cap \O Lu$, where $\bar u\in P/Q$ and $u$ is a representative of $\bar u$ in $P$.

\bigskip\noindent{\bf Lemma 5.9.}\quad {\it Keep the notation and the assumption of 5.8. Then $\tau$ is a $P/Q$-graded algebra isomorphism and $\pi$ is the restriction of $\tau$ to the $1$-components of $Z(\O Gb)$ and $Z(\O K)$.}

\medskip\noindent{\it Proof.}\quad
For any $r\in Z(\mathcal K Gb)$ and any $s\in Z(\K K)$, denote by $r(g)$ the class function on $G$ such that $r=\sum_{g\in G} r(g^{-1})g$, and by $s(t)$ the class function on $K$ such that $s=\sum_{t\in K}s(t^{-1})t$.
We define a $\mathcal K$-linear map $R_\nu^\circ: Z(\K G b)\rightarrow Z(\mathcal K K),\, r\mapsto \sum_{t\in K}\Big(\frac{1}{|G|}\sum_{g\in G}\nu(g,\,t)r(g)\Big)t^{-1}$, and a $\mathcal K$-linear
map
$I_\nu^\circ: Z(\mathcal K K)\rightarrow Z(\K G b),\,s\mapsto \sum_{g\in G}\Big(\frac{1}{|K|}\sum_{t\in K}\nu(g,\,t)s(t)\Big)g^{-1}$.
By \cite[Theorem 1.5]{B1}, the isomorphism $\tau$ maps $r$ onto $I_\nu^\circ (r R_\nu^\circ(b))$ for any $r\in Z(\O K)$. Since $\nu={\rm Ind}^{G\times K}_{(H\times L)\Delta(P)}(\nu')$, it is trivial to prove that $I_\nu^\circ$ and $R_\nu^\circ$ are $P/Q$-graded linear maps and thus that $\tau$ is a $P/Q$-graded algebra isomorphism.

We extend the isomorphisms $\tau$ and $\pi$ by the $\K$-linearity, and get $\K$-algebra isomorphisms $\tau^\K: Z(\K K)\cong Z(\K Gb)$ and $\pi^\K: Z(\K L)\cong Z(\K Hb)$. Note that $Z(\K L)$ is the 1-component of $Z(\K K)$. Since $\tau^K(Z(\K L))=Z(\K Gb)\cap Z(\K Hb)$ and
the algebras $Z(\K L)$ and $Z(\K Hb)$ have the same $\K$-dimension, we have $\tau^K(Z(\K L))=Z(\K Hb)$.

Since $P$ is abelian, by Lemma 5.2, any character $\chi\in {\rm Irr}(H, b)$ is extendible to $G$. Consequently,
for any $\chi\in {\rm Irr}(H, b)$ and any $\hat\chi\in {\rm Irr}(G, b)$, $\hat\chi$ is an extension of $\chi$ if and only if $e_\chi e_{\hat\chi}\neq 0$.

Fix $\chi\in {\rm Irr}(L)$ and $\hat\chi\in {\rm Irr}(K)$ such that $\hat\chi$ is an extension of $\chi$.
Since $e_\chi$ is primitive in $Z(\K L)$, $\tau^\K(e_\chi)=e_\varphi$ for a character $\varphi\in {\rm Irr}(H, b)$. Similarly, $\tau^\K(e_{\hat\chi})=e_{\hat\varphi}$ for some $\hat\varphi\in {\rm Irr}(G, b)$. By the proof of \cite[Theorem 1.5]{B1}, $I_\nu(\hat\chi)=\pm\hat\varphi$. Then we use the diagram 5.8.1 and obtain $I_{\nu'}(\chi)=\pm\varphi$. Again, by the proof of \cite[Theorem 1.5]{B1}, we have $\pi^\K(e_\chi)=e_\varphi$. The proof is done.

\bigskip\noindent{\bf Proposition 5.10.}\quad {\it Keep the notation and the assumption of 5.1. Assume that $E$ is cyclic and acts freely on $R-\{1\}$ and that there is a Morita equivalence between $\O Hb$ and $\O L$ compatible with the $\ast$-structure. Then there is a Morita equivalence between $\O Gb$ and $\O K$ compatible with the $\ast$-structure. }

\medskip\noindent{\it Proof.}\quad
We suppose that the $\O(H\times L)$-module $M$ induces the Morita equivalence between $\O Hb$ and $\O L$, and denote by $\mu'$ the character of $M$. Then the map $I_{\mu'}: \mathbb Z{\rm Irr}(L)\rightarrow \mathbb Z{\rm Irr}(H, b)$ is a perfect isometry compatible with the $\ast$-structure, which induces an algebra isomorphism $$\sigma: Z(\O L)\cong Z(\O Hb).$$
By Proposition 5.7,
the character $\mu'$ has a suitable extension $\tilde\mu'$ to $(H\times L)\Delta(P)$ such that setting $\mu={\rm Ind}^{G\times K}_{(H\times L)\Delta(P)}(\tilde\mu')$, the map $I_\mu: \mathbb Z{\rm Irr}(K)\rightarrow \mathbb Z{\rm Irr}(G, b)$ is a perfect isometry compatible with the $\ast$-structure and
the diagram $$\begin{array}{ccc}
  \mathbb Z{\rm Irr}(K) & \xrightarrow[I_\mu]{}  &  \mathbb Z{\rm Irr}(G, b) \\
  \downarrow &  & \downarrow \\
  \mathbb Z{\rm Irr}(L) & \xrightarrow{I_{\mu'}}  & \mathbb Z{\rm Irr}(H, b).
\end{array}$$ is commutative, where the two down arrows are induced by the restrictions of characters.
Since $I_{\mu'}$ maps ordinary characters of $L$ onto ordinary characters of $H$, by the diagram, $I_\mu$ maps ordinary characters of $K$ onto ordinary characters of $G$. In particular, $\mu$ is an ordinary character of $G\times K$.

The map $I_\mu$ induces an algebra isomorphism $$\varsigma: Z(\O K)\cong Z(\O Gb).$$
By Lemma 5.9, the isomorphism $\varsigma$ is a $P/Q$-graded algebra isomorphism and the restriction of $\varsigma$ to the 1-component of $Z(\O Gb)$ coincides with $\sigma$. Since $T$ is contained in $Z(\O K)$, for any $t\in T$, $\varsigma(t)$ is contained in the $\bar t$-component of $Z(\O Gb)$ and $t\zeta(t^{-1})$ is inside $\O Hb$. The map $\rho: T\rightarrow (\O Hb)^*, t\mapsto t\varsigma(t^{-1})$ is a group homomorphism. Now we extend the $\O(H\times L)$-module $M$ to $(H\times L)\Delta(P)$ by the equality $(t, t)\cdot m=\rho(t)m$ for any $t\in T$ and any $m\in M$. We claim that such an extension is well defined.
Since $\O(H\times L)$-module $M$ induces the Morita equivalence between $\O Hb$ and $\O L$, by Lemma 2.10, we have $\sigma(a)m=ma$ for any $a\in Z(\O L)$ and any $m\in M$. So for any $t\in T\cap H$ and any $m\in M$, we have $$(t, t)m=tmt^{-1}=t\sigma(t^{-1})m=t\varsigma(t^{-1})m=(t, t)\cdot m.$$ The claim is done. Now by \cite[Theorem 3.4 (a)]{MA}, the $\O(G\times K)$-module ${\rm Ind}^{G\times K}_{(H\times L)\Delta(P)}(M)$ induces a Morita equivalence between $\O Gb$ and $\O K$.
Now it remains to show that the Morita equivalence is compatible with the $\ast$-structure.

Set $\tilde M={\rm Ind}^{G\times K}_{(H\times L)\Delta(P)}(M)$. The Morita equivalence induced by $\tilde M$ induces an algebra isomorphism (see Paragraph 2.9) $$\rho:Z(\O K)\cong Z(\O Gb)$$ such that for any $a\in Z(\O Gb)$ and any $a'\in Z(\O K)$, $a$ and $a'$ correspond to each other if and only if $am=ma'$ for any $m\in \tilde M$. We claim that $\rho$ and $\varsigma$ are the same.
For any $a\in Z(\O L)$ and any $m\in M$, we have $$(\sigma(a)\otimes 1)\otimes m=(1\otimes 1)\otimes \sigma(a)m=(1\otimes 1)\otimes ma=(1\otimes a^\circ)\otimes m,$$ where $a^\circ$ is the image of $a$ through the opposite ring isomorphism $\O L\rightarrow \O L$ sending $x$ onto $x^{-1}$ for any $x\in L$.  So $\rho(a)=\sigma(a)=\varsigma(a)$.
On the other hand, given $t\in T$ and $m\in M$, we have $$(\varsigma(t)\otimes t)\otimes m=\Big((t\otimes t)(t^{-1}\varsigma(t)\otimes 1)\Big)\otimes m=(t\otimes t)\otimes t^{-1}\varsigma(t)m=(1\otimes 1)\otimes m $$ and thus $(\varsigma(t)\otimes 1)\otimes m=(1\otimes t^{-1})\otimes m$. So $\rho$ maps $t$ onto $\varsigma(t)$. Since $Z(\O K)=\sum_{t\in T}Z(\O L)t$, the claim is done.

Finally, we claim that the character of ${\rm Ind}^{G\times K}_{(H\times L)\Delta(P)}(M)$ is equal to $\mu$. In particular, the Morita equivalence between $\O Gb$ and $\O K$ induced by ${\rm Ind}^{G\times K}_{(H\times L)\Delta(P)}(M)$ induces a perfect isometry compatible with the $\ast$-structure.
Indeed, since the isomorphisms $\rho$ and $\varsigma$ coincide with each other, setting $\rho(e_\chi)=e_{\chi'}$ for any $\chi\in {\rm Irr}(K)$, by the proof of \cite[Theorem 1.5]{B1},
both $\mu$ and the character of $\tilde M$ are equal to the sum $\sum_{\chi\in {\rm Irr}(K)} \chi'\times \chi^\ast$. The claim is done.

\bigskip\bigskip\noindent{\bf 6. Proof of the Theorem}

\bigskip In this section, we prove {\bf Theorem} in the introduction. We borrow the notation and the assumption in {\bf Theorem}. Note that $p$ is 2 and that the defect group $P$ is abelian. However, some results, such as Lemmas 6.2, 6.5 and 6.6, hold in general. Let $(P, f)$ be a maximal Brauer pair associated with the block $b$. The inertial quotient $E_G(P, f)$ may be of order 1, 3, 7 or 21.

\bigskip\noindent{\bf Lemma 6.1.}\quad {\it Let $n$ be a nilpotent block of a normal subgroup $N$ of the group $G$ covered by the block $b$, and set $N'=NZ(G)$. Then there is a nilpotent block $n'$ of $N'$ covered by the block $b$. }

\medskip\noindent{\it Proof.}\quad We consider the group homomorphism $N\times Z(G)\rightarrow N', (x, y)\mapsto xy$ and the induced algebra homomorphism $\O(N\times Z(G))\rightarrow \O N'$. Since the homomorphism $\O(N\times Z(G))\rightarrow \O N'$ maps $n\otimes 1$ onto $n$, by the same proof method as that of Lemma 2.1, we prove that any block of $N'$ covering $n$ is the image of the block $n\otimes w$ of $N\times Z(G)$ for some suitable block $w$ of $Z(G)$, that its inertial quotient has order 1, and that it is nilpotent.

\bigskip\noindent{\bf Lemma 6.2.}\quad {\it Assume that there is a finite group $\tilde G$, such that $G$ is normal in $\tilde G$ and $\tilde G$ has a nilpotent block $\tilde b$ covering $b$. Then $\O Gb$ and $\O N_G(P)c$ are basically Morita equivalent. }

\medskip\noindent{\it Proof.}\quad By the main theorem of \cite{P1}, the block $b$ is inertial. That is to say, $\O Gb$ and $\O N_G(P)c$ are basically Morita equivalent.

\bigskip\noindent{\bf Lemma 6.3.}\quad {\it Let $G$ be a finite group and $N$ a normal subgroup of $G$ with index an odd prime. Let $b$ be a block of $G$ with defect group $P$, and $\frak n$ a $G$-stable block of $N$ covered by $b$. Assume that $P$ is $\mathbb Z_{2^n}\times \mathbb Z_{2^n}\times \mathbb Z_{2^n}$, that the inertial quotient of the block $b$ is of order $7$ or $21$ and that the block $\frak n$ is not nilpotent. Then the inertial quotient of the block $\frak n$ is of order $7$ or $21$. }

\medskip\noindent{\it Proof.}\quad Since $N$ is a normal subgroup of $G$ with index an odd prime and $b$ covers $\frak n$, $P$ is also a defect group of $\frak n$; moreover, we can choose a maximal $\frak n$-Brauer pair $(P, h)$ so that $f$ covers $h$.

Suppose that $N$ contains $C_G(P)$. Then $f=h$, $N_N(P, h)/C_N(P)$ is a normal subgroup of $N_G(P, f)/C_G(P)$, and $N_N(P, h)/C_N(P)$ has order 7 or 21 since the block $\frak n$ is not nilpotent.

Suppose that $N$ does not contain $C_G(P)$ and that $C_G(P)=C_G(P, h)$, where $C_G(P, h)$ is the intersection of $C_G(P)$ and the stabilizer $N_G(P, h)$ of $(P, h)$ under the $G$-conjugation. Then $h$ is a central idempotent in $\O C_G(P)$ and we have $hf=f$.
For any $x\in N_G(P, f)$, we have $xhx^{-1}f=f$, $xhx^{-1}=h$ and $x\in N_G(P, h)$. So $N_G(P, f)\leq N_G(P, h)$. Since the index of $N$ in $G$ is an odd prime and $N$ does not contain $C_G(P)$, $G=N C_G(P)$ and $N_G(P, h)=N_N(P, h)C_G(P)$. The inclusion $N_G(P, f)\subset N_G(P, h)$ induces an injective group homomorphism $$N_G(P, f)/C_G(P)\rightarrow N_G(P, h)/C_G(P)\cong N_N(P, h)/C_N(P).$$ This implies that the inertial quotient of the block $\frak n$ has order 7 or 21.

Suppose that $N$ does not contain $C_G(P)$ and that $C_G(P)\neq C_G(P, h)$. Since the index of $N$ in $G$ is an odd prime $q$, the index of $C_N(P)$ in $C_G(P)$ is $q$ and $C_N(P)=C_G(P, h)$. That is to say, the stabilizer of $h$ under the $C_G(P)$-conjugation is $C_N(P)$. Since $f$ covers $h$, we have $f=\sum_{x\in I} xhx^{-1}$, where $I$ is a complete set of representatives of $C_N(P)$ in $C_G(P)$ and $xhx^{-1} h=0$ for any $x\in I$ outside $C_N(P)$. Obviously we have $N_G(P, h)\subset N_G(P, f)$. For any $y\in N_G(P, f)$, there is $z\in C_G(P)$ such that $z^{-1}y$ centralizes $h$. So $N_G(P, f)\subset N_G(P, h)C_G(P)$ and then $N_G(P, f)= N_G(P, h)C_G(P)$. Now it is easily seen that the inclusion $N_N(P, h)\subset N_G(P, f)$ induces an injective group homomorphism $$N_N(P, h)/C_N(P)\rightarrow N_G(P, f)/C_G(P),$$ whose image is normal. So the inertial quotient of the block $\frak n$ has order 7 or 21.

\bigskip\noindent{\bf Proposition 6.4.}\quad {\it Assume that $P$ is $\mathbb Z_{2^n}\times \mathbb Z_{2^n}\times \mathbb Z_{2^n}$ and that $E_G(P, f)$ is of order $7$ or $21$. If $n\geq 2$, then $\O Gb$ and $\O N_G(P)c$ are basically Morita equivalent.}

\medskip\noindent{\it Proof.}\quad Suppose that $b$ is a block of $G$ with defect group $P$ and the inertial quotient of order 7 or 21,
such that $(|G:Z(G)|, |G|)$ is minimal in the lexicographic ordering and that $\O Gb$ and $\O N_G(P)c$ are not basically Morita equivalent. Then the block $b$ has to be quasiprimitive.

Suppose that $O_2(G)$ is not trivial. Set $H=C_G(O_2(G))$. Let $h$ be a block of $H$ covered by the block $b$. Then $P$ is a defect group of the block $h$. Since $C_G(P)\subset H$, we may adjust the choice of $h$ so that $(P, f)$ is a maximal Brauer pair of the block $h$. By \cite[Proposition 15.10]{P2}, $h$ is equal to $b$ and the quotient group $G/H$ is an odd group. The inertial quotient $E_H(P, f)$ is a normal subgroup of $E_G(P, f)$. Since $E_G(P, f)$ is of order 7, $E_H(P, f)$ is of order 1, or 7, or 21. Suppose that $E_H(P, f)$ is of order 1. Then the block $b$ of $H$ is nilpotent and by \cite[Theorem]{Z2}, $\O Gb$ and $\O N_G(P)c$ are basically Morita equivalent. That contradicts the choice of the block $b$. Suppose that $E_H(P, h)$ is of order 7 or 21. Then the commutator subgroup of $P$ and $E_H(P, h)$ is $P$. But since $O_2(G)$ is in the center of $H$, that is impossible.
So $O_2(G)$ has to be trivial and $Z(G)$ is a subgroup of odd order.

Suppose that the block $b$ covers a nilpotent block $n$ of a normal subgroup $N$ of $G$. Set $N'=N Z(G)$. By Lemma 6.1, there is a nilpotent block $n'$ covered by the block $b$. By Lemma 4.5 applied to the normal subgroup $N'$ and its block $n'$ and by the minimality of $(|G:Z(G)|, |G|)$,
$N'$ has to be equal to $Z(G)$ and thus $N$ is contained in $Z(G)$.

Denote by $E(G)$ the layer of the group $G$ and by $F(G)$ the Fitting subgroup of $G$. By the choice of the block $b$, $F(G)$ is equal to $O_{2'}(G)$. Set $F^*(G)=E(G) F(G)$, the generalized Fitting subgroup of $G$. Suppose that $E(G)$ is the central product $L_1*L_2*\cdots *L_t$, where $L_1, L_2,\cdots, L_t$ are components of $G$. Since $C_G(F^*(G))\leq F^*(G)$, by the choice of the block $b$, $t\geq 1$.

Let $e$ be the unique block of $E(G)$ covered by the block $b$ and $e_i$ a block of $L_i$ covered by the block $e$. Suppose $t>1$.
Since the block $b$ has defect group $\mathbb Z_{2^n}\times \mathbb Z_{2^n}\times \mathbb Z_{2^n}$, by Lemma 2.1 some block $e_i$ has to be nilpotent and then the block $b$ covers a nilpotent block of a normal subgroup generated by all $G$-conjugates of $L_i$. So $L_i\leq Z(G)$ and this contradicts $L_i$ being a component of $G$. So $t=1$ and then $F^*(G)=L_1O_{2'}(G)$.

By the Schreier conjecture, $G/L_1$ is solvable. We claim that $G/L_1$ is trivial.

Suppose that $G$ has a normal subgroup $H$ of index 2. By \cite[Proposition 5.3]{KP} $G=PH$. Since $P=[P, E_G(P, f)]$,
by \cite[Proposition 4.2]{P6} $P$ is contained in $H$ and so $G=H$. That contradicts the choice of $H$. So $G$ has no normal subgroup $H$ of index 2.

Suppose that $G$ has a normal subgroup $H$ with an odd prime index containing $L_1$. Let $h$ be the unique block of $H$ covered by $b$. Note that $P$ is a defect group of the block $h$.
By \cite[Theorem]{Z2}, the block algebra $\O Hh$ is not basically Morita equivalent to its Brauer correspondent in $N_H(P)$. By Lemma 6.3, the inertial quotient of the block $h$ has order 7 or 21. This is against the minimality of $(|G:Z(G)|, |G|)$.

So $G$ is equal to $L_1$. Since $P$ is $\mathbb Z_{2^n}\times \mathbb Z_{2^n}\times \mathbb Z_{2^n}$ for $n\geq 2$,
only the case (iii) among the four cases of \cite[Theorem 6.1]{EKKS} happens for the block $b$ of $G$. In this case, there is a finite group $\tilde G$, such that $G$ is normal in $\tilde G$ and $\tilde G$ has a nilpotent block $\tilde b$ covering $b$.
By Lemma 6.2, $\O Gb$ and $\O N_G(P)c$ are basically Morita equivalent. That contradicts the choice of the block $b$.

\bigskip\noindent{\bf Lemma 6.5.}\quad {\it Let $H$ and $H'$ be finite groups and $h$ and $h'$ blocks of $H$ and $H'$. Assume that $M$ is an indecomposable $\O(H\times H')$-module inducing a Morita stable equivalence between $\O Hh$ and $\O H'h'$. Let $\ddot R$ be a vertex of $M$ and $R$ and $R'$ the images of $\ddot R$ through the projections $H\times H'\rightarrow H$ and $H\times H'\rightarrow H'$. Then $R$ and $R'$ are defect groups of $h$ and $h'$. }

\medskip\noindent{\it Proof.}\quad This lemma follows from  \cite[Theorem 6.9]{P3}.

\bigskip\noindent{\bf Lemma 6.6.}\quad {\it Let $H$, $H'$ and $H''$ be finite groups, and $h$, $h'$ and $h''$ blocks of $H$, $H'$ and $H''$. Assume that $R$ is a common $p$-subgroup of $H$, $H'$ and $H''$, that $M$ is an indecomposable $\O(H\times H')$-module with vertex $\Delta_\sigma(R)=\{(u, \sigma(u))|u\in R\}$ for some group automorphism $\sigma$ on $R$, and that $M'$ is an indecomposable $\O(H'\times H'')$-module with vertex $\Delta_{\sigma'}(R)=\{(u, \sigma'(u))|u\in R\}$ for some group automorphism $\sigma'$ on $R$. Then the order of vertex of any indecomposable direct summand of $M\otimes_{\O H'} M'$ is at most $|R|$. }

\medskip\noindent{\it Proof.}\quad Let $\O \Delta_\sigma(R)$-module $S$ be a source of $M$ and $\O \Delta_{\sigma'}(R)$-module $S'$ a source of $M'$. Then $M$ is a direct summand of ${\rm Ind}^{H\times H'}_{\Delta_\sigma(R)}(x)$ and $M'$ is a direct summand of ${\rm Ind}^{H'\times H''}_{\Delta_{\sigma'}(R)}(S')$.
We have an $\O(H\times H'')$-module isomorphism
$${\rm Ind}^{H\times H''}_{R\times R}(\O H'\otimes_\O {\rm Res}_\tau(S')\otimes_\O {\rm Res}_{\tau'}(S''))\cong {\rm Ind}^{H\times H'}_{\Delta_\sigma(R)}(x)\otimes_{\O H'} {\rm Ind}^{H'\times H''}_{\Delta_{\sigma'}(R)}(S'), $$ mapping $( x\otimes y)\otimes (z\otimes s\otimes s')$ onto $((x\otimes z)\otimes s')\otimes((1\otimes y)\otimes s')$ for any $x\in \O H$, any $y\in \O H''$, any $z\in \O H'$, any $s\in S$ and $s'\in S'$, where $\tau$ is the homomorphism $R\times R\rightarrow \Delta_\sigma(R), (u, v)\mapsto (u, \sigma(u))$, $\tau'$ is the homomorphism $R\times R\rightarrow \Delta_{\sigma'}(R), (u, v)\rightarrow (\sigma'^{-1}(v), v)$ and $\O H'$ is an $\O(R\times R)$-module defined by the equality $$(u, v)w=\sigma(u)w\sigma'(v^{-1})$$ for any $u, v\in R$ and any $w\in \O H$. Since $\O H'\cong \oplus_{x\in R\backslash H'/R} {\rm Ind}^{R\times R}_{R_x}(\O)$, where $R\backslash H'/R$ denotes a set of representatives of the double cosets of $R$ in $H'$ and $R_x=\{(\sigma^{-1}(z), \sigma'^{-1}(x^{-1}zx))|z\in xRx^{-1}\cap R\}$ for any $x\in R\backslash H'/R$, the order of vertex of any indecomposable direct summand of ${\rm Ind}^{H\times H'}_{\Delta_\sigma(R)}(x)\otimes_{\O H'} {\rm Ind}^{H'\times H''}_{\Delta_{\sigma'}(R)}(S')$ is at most $|R|$, and so is that of $M\otimes_{\O H'} M'$.

\bigskip\noindent{\bf Proposition 6.7.}\quad {\it Let $h$ be a block of $H$ with defect group $R=\mathbb Z_{2^n}\times \mathbb Z_{2^n}$, where $n\geq 2$. Then the block algebra $\O Hh$ and its Brauer correspondent $\O N_H(R)d$ are basically Morita equivalent.}

\medskip\noindent{\it Proof.}\quad
Let $(R, e)$ be a maximal $h$-Brauer pair and set $K=R\rtimes E_H(R, e)$. By the structure theorem of blocks with normal defect groups, there is a Morita equivalence between $\O N_H(R)d$ and $\O K$ induced by a $p$-permutation bimodule.

By \cite[Theorem 1.1]{EKKS}, $\O Hh$ and $\O K$ are Morita equivalent. Suppose that an $\O(H\times K)$-module $M$ induces the Morita equivalence. By \cite[Remarque 6.8]{P9} there is a stable equivalence between $\O Hh$ and $\O K$ induced by an indecomposable $\O(H\times K)$-module $\tilde M$ with endopermutation module sources and vertex $\Delta(R)$. Denote by $\tilde M^*$ the dual of the $\O(H\times K)$-module $\tilde M$. Then $\tilde M^*\otimes_{\O G} M$ is an indecomposable $\O(K\times K)$-module inducing a Morita stable self-equivalence on $\O K$. By \cite[Corollary 3.3]{CR}, there is an integer $r$ such that $\Omega^r(\tilde M^*\otimes_{\O G} M)$ induces a Morita self-equivalence on $\O K$, where $\Omega$ denotes the Heller translate.
By \cite[Theorem 2]{Z3} and \cite[Corollary 7.4]{P3}, $\Omega^r(\tilde M^*\otimes_{\O G} M)$ has vertex $\Delta_\sigma(R)=\{(u, \sigma(u))|u\in R\}$ for some group automorphism $\sigma$ on $R$, and so does $\tilde M^*\otimes_{\O G} M$. Since $M$ is a direct summand of $\tilde M\otimes_{\O K} \tilde M^*\otimes_{\O G} M$ and the order of vertex of any indecomposable direct summand of $\tilde M\otimes_{\O K} \tilde M^*\otimes_{\O G} M$ is at most $|R|$ (see Lemma 6.6), by Lemma 6.5 the order of vertex of $M$ has to be equal to $|R|$. Then by \cite[Corollary 7.4]{P3}, $M$ induces a basic Morita equivalence between $\O Hh$ and $\O K$.

\bigskip\noindent{\bf Lemma 6.8.}\quad {\it Let $G$ be a finite group and $N$ a normal subgroup of $G$ with index an odd prime. Let $b$ be a block of $G$ with defect group $P$, and $\frak n$ a $G$-stable block of $N$ covered by $b$. Assume that $P$ is $\mathbb Z_{2^n}\times \mathbb Z_{2^n}\times \mathbb Z_{2^m}$, that the inertial quotient of the block $b$ is of order $3$ and that the block $\frak n$ is not nilpotent. Then the inertial quotient of the block $\frak n$ has order $3$ or $21$.}

\medskip\noindent{\it Proof.}\quad One uses the proof of Lemma 6.3 to prove the lemma.

\bigskip\noindent{\bf Lemma 6.9.}\quad {\it Let $G$ be a finite group and $b$ a block of $G$ with defect group $P$. Assume that $P$ is $\mathbb Z_{2^n}\times \mathbb Z_{2^n}\times \mathbb Z_{2^m}$ and that the inertial quotient of the block $b$ is of order $3$. Then $\ell_G(b)$ is $3$.}

\medskip\noindent{\it Proof.}\quad See \cite[Theorem 1.1 (i)]{HZ} or alternatively use \cite[Theorem 1]{W1} to prove the lemma.

\bigskip\noindent{\bf Proposition 6.10.}\quad {\it Assume that $P$ is $\mathbb Z_{2^n}\times \mathbb Z_{2^n}\times \mathbb Z_{2^m}$ and that $E_G(P, f)$ is of order $3$.
If $n\geq 2$, then there is a Morita equivalence between the block algebras $\O Gb$ and $\O N_G(P)c$ compatible with the $\ast$-structure, where $c$ denotes the Brauer correspondent of $b$ in $N_G(P)$.}

\medskip\noindent{\it Proof.}\quad Set $K=P\rtimes E_G(P, f)$. By the structure theorem of blocks with normal defect groups, there is a Morita equivalence between $\O N_G(P)c$ and $\O K$ induced by a $p$-permutation bimodule. Therefore, in order to prove the lemma, by Propositions 2.7 and 3.6, it suffices to show that there is a Morita equivalence between $\O Gb$ and $\O K$ compatible with the $\ast$-structure.

When $m=0$, by Propositions 6.7, 2.7 and 3.6 there is a Morita equivalence between $\O Gb$ and $\O K$ compatible with the $\ast$-structure. We go by induction on $m$. Suppose that $b$ is a block of $G$ with defect group $P$ and the inertial quotient of order 3
such that $(|G:Z(G)|, |G|)$ is minimal in the lexicographical order and that there is not a Morita equivalence between $\O Gb$ and $\O K$ compatible with the $\ast$-structure. By Lemma 3.4, the block $b$ has to be quasiprimitive.

Set $H=C_G(O_2(G))$. Let $h$ be the block of $H$ covered by the block $b$. Then $P$ is a defectgroup of the block $h$. Since $C_G(P)\subset H$, we may assume without loss that $(P, f)$ is a maximal Brauer pair of the block $h$. By \cite[Proposition 15.10]{P2}, $h$ is equal to $b$ and the quotient group $G/H$ is an odd group. The inertial quotient $E_H(P, f)$ is a normal subgroup of $E_G(P, f)$. Since $E_G(P, f)$ is of order 3, $E_H(P, f)$ is of order 1 or 3. If $E_H(P, f)$ is of order 1, then the block $b$ of $H$ is nilpotent and by \cite{Z2}, $\O Gb$ and $\O K$ are basically Morita equivalent. By Proposition 3.6, there is a Morita equivalence between $\O Gb$ and $\O K$ compatible with the $\ast$-structure. That contradicts the choice of the block $b$. So $E_H(P, f)$ is of order 3 and $N_G(P, f)=N_H(P, f)C_G(P)$. By a Frattini argument, we have $G=N_G(P, f)H$. So $G=H$ and $O_2(G)$ is central.

Suppose that the block $b$ covers a nilpotent block of a normal subgroup $N$ of $G$. Set $N'=N Z(G)$. By Lemma 6.1, there is a nilpotent block $n'$ covered by the block $b$. By Lemma 3.5 applied to the normal subgroup $N'$ and its block $n'$ and the minimality of $(|G:Z(G)|, |G|)$,
$N'$ has to be equal to $Z(G)$ and thus $N$ is contained in $Z(G)$.

Denote by $E(G)$ the layer of the group $G$ and by $F(G)$ the Fitting subgroup of $G$. Then $F(G)$ is equal to $Z(G)$. Set $F^*(G)=E(G) F(G)$, the generalized Fitting subgroup of $G$. Suppose that $E(G)$ is the central product $L_1*L_2*\cdots *L_t$, where $L_1, L_2,\cdots, L_t$ are the components of $G$. Since $C_G(F^*(G))\leq F^*(G)$, by the choice of the block $b$, $t\geq 1$.

Let $e$ be the unique block of $E(G)$ covered by the block $b$ and $e_i$ a block of $L_i$ covered by the block $e$.
By the choice of the block $b$, the order of the inertial quotient of the block $e$ may be 3, 7, or 21.
Suppose $t>1$.
By Lemma 2.1, some block $e_i$ has to be nilpotent and then the block $b$ covers a nilpotent block of a normal subgroup generated by all $G$-conjugates of $L_i$. So $L_i$ is contained in $Z(G)$. This contradicts $L_i$ being a component of $G$. So $t=1$ and $F^*(G)=L_1Z(G)$.

By the Schreier conjecture, $G/L_1$ is solvable. We claim that $G/L_1$ is trivial.

Suppose that $G$ has a normal subgroup $H$ of index 2. Since $E_G(P, f)$ is of order 3, we may assume without loss of generality that the commutator $[P, E_G(P, f)]$ is equal to $\mathbb Z_{2^n}\times \mathbb Z_{2^n}$ and that the subgroup of $E_G(P, f)$-fixed elements of $P$ is $\mathbb Z_{2^m}$. Let $h$ be the unique block of $H$ covered by the block $b$. Then the intersection $Q=P\cap H$ is a defect group of the block $h$. By \cite[Proposition 4.2]{P6} $\mathbb Z_{2^n}\times \mathbb Z_{2^n}\leq Q$ and so $Q$ is $\mathbb Z_{2^n}\times \mathbb Z_{2^n}\times \mathbb Z_{2^{m-1}}$ since $G=PH$ (see \cite[Proposition 5.3]{KP}). By \cite[Lemma 3.6]{Z1} the block $h$ has the inertial quotient of order 3. By induction on $m$, there is a Morita equivalence between $\O Hh$ and $\O L$ compatible with the $\ast$-structure, where $L=Q\rtimes E_G(P, f)$. By Proposition 5.10 there is a Morita equivalence between $\O Gb$ and $\O K$ compatible with the $\ast$-structure. This contradicts the choice of the block $b$. So $G$ has no normal subgroup $H$ of index 2.

Suppose that $G$ has a normal subgroup $H$ of an odd prime index containing $L_1$. Let $h$ be the unique block of $H$ covered by the block $b$.
By Lemma 6.8, the inertial quotient of the block $h$ has order 3 or 21. By Proposition 6.4 and \cite[Theorem]{Z2}, we may exclude the case that the inertial quotient of the block $h$ has order 21. Now by Lemma 6.9, $\ell_G(b)=\ell_H(h)=3$.
Since $H$ has odd prime index in $G$, by the last paragraph of the proof of \cite[Theorem 1.1]{CEKL}, that is impossible.

So $G=L_1$. Then as in the last paragraph of the proof of Proposition 6.4, we prove that the block algebra $\O Gb$ and $\O K$ are basically Morita equivalent. Then by Propositions 2.7 and 3.6, there is a Morita equivalence between $\O Gb$ and $\O K$ compatible with the $\ast$-structure. That contradicts the choice of the block $b$.

\bigskip\noindent{\bf 6.11.}\quad Proof of {\bf Theorem}

\medskip Note that $\O$ unnecessarily has characteristic 0 in {\bf Theorem}. We divide the proof of {\bf Theorem} into the characteristic zero case and characteristic nonzero case.

Firstly, we assume that $\O$ has characteristic 0. When $E_G(P, f)$ is 1, the block $b$ of $G$ is nilpotent and {\bf Theorem} is true (see \cite{P4}). When $E_G(P, f)$ is 7 or 21, $n=m$ and {\bf Theorem} follows from Proposition 6.4. When $E_G(P, f)$ is 3, {\bf Theorem} follows from Proposition 6.10.

Finally, we assume that $\O$ has characteristic $p$. The blocks $b$ and $c$ determine blocks $\bar b$ and $\bar c$ of $kG$ and $kN_G(P)$ with defect group $P$. By \cite[Chapter II, Theorem 3]{S}, there is a complete discrete valuation ring $\tilde \O$ with characteristic 0 and with residue field $k$. The blocks $\bar b$ and $\bar c$ can be lifted to blocks $\tilde b$ and $\tilde c$ of $\tilde\O G$ and $\tilde \O N_G(P)$ with defect group $P$. By \cite[Chapter II, Proposition 3]{S}, without loss of generality, we may assume that $\tilde \O$ contains a $|G|$-th primitive root of unity.
By the last paragraph, the block algebras $\tilde\O G\tilde b$ and $\tilde\O N_G(P)\tilde c$ are Morita equivalent. So the block algebras $k G \bar b$ and $k N_G(P)\bar c$ are Morita equivalent. Since $\O$ and $k$ have the same characteristic $p$, by \cite[Chapter II, Proposition 8]{S}, $k$ can be identified with a subring of $\O$. Moreover, it is easy to see that $\O Gb$ and $\O N_G(P)c$ are equal to $\O\otimes_k k G\bar b$ and $\O\otimes_k k N_G(P)\bar c$ respectively. Therefore the block algebras $\O Gb$ and $\O N_G(P)c$ are Morita equivalent.



\end{document}